\documentclass{amsart}

\usepackage{kotex}
\usepackage{amsthm}
\usepackage{amssymb}
\usepackage{amsmath}
\usepackage{amscd}
\usepackage{enumerate}
\usepackage{graphicx} 
\usepackage{siunitx}
\usepackage{tikz-cd}
\usepackage{bm}
\usepackage{dsfont}

\numberwithin{equation}{section}

\usepackage{verbatim}

\usepackage{hyperref}
\hypersetup{
    colorlinks=true,
    linkcolor=blue,
    filecolor=magenta,      
    urlcolor=cyan,
    citecolor=red}

\usepackage{mathtools}
\usepackage[tableposition=top]{caption}
\usepackage{booktabs,dcolumn}

\newtheorem{theorem}{Theorem}[section]
\newtheorem{proposition}[theorem]{Proposition}
\newtheorem{lemma}[theorem]{Lemma}
\newtheorem{corollary}[theorem]{Corollary}
\newtheorem{conjecture}[theorem]{Conjecture}

\newtheorem{definition}[theorem]{Definition}
\newtheorem{remark}[theorem]{Remark}

\newcommand{\bB}{\mathbb{B}}

\newcommand{\bR}{\mathbb{R}}
\newcommand{\bZ}{\mathbb{Z}}
\newcommand{\bN}{\mathbb{N}}

\newcommand{\cF}{\mathcal{F}}

\newcommand{\vl}{{\vec\ell}}

\newcommand{\supp}{\mathrm{supp}}
\newcommand{\dist}{\mathrm{dist}}

\begin{document}

\title{Maximal averages and non-transversality}

\author[J.~B. Lee]{Jin Bong Lee}
\address{Research Institute of Mathematics, Seoul National University, 08826 Gwanak-ro 1, Seoul, Republic of Korea}
\email{jinblee@snu.ac.kr}

\author[J. Lee]{Juyoung Lee}
\address{Department of Mathematics, Korea Institute for Advanced Study, Seoul 02455, Republic of Korea}
\email{juyounglee@kias.re.kr}

\author[J. Oh]{Jeongtae Oh}
\address{Research Institute of Mathematics, Seoul National University, 08826 Gwanak-ro 1, Seoul, Republic of Korea}
\email{ojt0117@snu.ac.kr}

\author[S. Oh]{Sewook Oh}
\address{June E Huh Center for Mathematical Challenges, Korea Institute for Advanced Study, Seoul 02455, Republic of Korea}
\email{sewookoh@kias.re.kr}


\begin{abstract}

We investigate the $L^p$ mapping properties of maximal functions associated with analytic hypersurfaces in $\mathbb R^d$, with a particular emphasis on the role of transversality.
Around points that are not transversal, we show that the associated maximal function is bounded on $L^p(\mathbb R^d)$ for all $p>2$, regardless of the decay of the Fourier transform of surface measures. In contrast, away from \emph{non-transversal} points, we prove that $L^p$ bounds for the maximal operator imply that the Fourier transform of the surface measure decays at rate $1/q$ for  $q>p$. Combining these two regimes, we demonstrate that the conjecture of Stein and Iosevich--Sawyer on maximal functions could be re-formulated, in the analytic setting, by restricting attention to transversal points. Moreover, our result completely settles the refined form of the conjecture for certain cases. 
\end{abstract}

\maketitle

\section{Introduction}

Let $\Gamma$ be a smooth compact hypersurface in $\bR^{d}$ for $d\ge2$.
For the induced Lebesgue measure $\mu$ on $\Gamma$, we define the isotropic dilation $\mu_t$ of the measure $\mu$ by
\[ \langle \mu_t, f \rangle =\int_\Gamma f(t x)\,\mathrm d \mu(x) . \]
Then, consider the associated maximal function given by
\begin{align}
\nonumber
    M_\Gamma f(x)\coloneq\sup_{t>0}|A_\Gamma f(t, x)| \coloneq  \sup_{t>0} |\mu_t \ast f(x)|.
\end{align}
Studies on maximal averages, $M_\Gamma$, have been extensively pursued for the last decades. 
Particularly, many works have been devoted to the case of $\Gamma$ being a hypersurface. 
The purpose of this work is to clarify the connection between $L^p$ estimates for maximal functions and geometric conditions of surfaces under which such estimates hold. By one of such conditions, we mean a decay of the Fourier transform of surface-carried measures. Moreover, we denote a Fourier decay of $\mu$ of rate $\rho$ as 
\begin{align}\label{ineq_fdecay}
    |\widehat{\mathrm d \mu}(\xi)| \le C (1+|\xi|)^{-\rho},\quad\text{for some $C>0$}.
\end{align}

It was E.~M. Stein \cite{St1976} who initiated the systematic study of maximal functions associated with submanifolds, beginning with the case of maximal averages over spheres in dimensions $d \ge 3$.
Greenleaf \cite{Gr1981} proved that if a hypersurface has at least $k$ nonvanishing principal curvatures, then Stein's $L^2$ method based on the Fourier decay of surface measures gives analogous results for $k\ge 2$. 
Later, it was further generalized by Rubio de Francia \cite{RdF1986} to measures that satisfy \eqref{ineq_fdecay} with the rate $\rho>1/2$.
For $\Gamma = \mathbb S^1$ in $\mathbb R^2$, \eqref{ineq_fdecay} holds  with $\rho = 1/2$. It is proved by Bourgain \cite{Bo1986} that the circular maximal function is bounded on $L^p(\mathbb R^2)$ for $p>2$, thereby establishing the sharp range.

For general hypersurfaces, however, the problem of determining the sharp maximal bounds remains widely open. The most complete picture so far arises in the convex case. In this context, Nagel--Seeger--Wainger \cite{NSW}, following the pioneering work of Bruna--Nagel--Wainger \cite{BNW}, developed a detailed correspondence between the $L^p$ boundedness of the maximal operator and the geometric features of the surface. Later, Iosevich--Sawyer \cite{IS2} confirmed this relationship for the case when the decay inequality \eqref{ineq_fdecay} holds with the rate $\rho \le 1/2$. 
Moreover, the authors \cite{IS2} suggested the following conjecture, which is now known as the Iosevich--Sawyer--Stein conjecture.

\begin{conjecture}[{\cite[ Conjecture~2]{IS2}}, \cite{IKM}]\label{conj_stein}
    Suppose $\Gamma$ is a smooth hypersurface and a measure $\mu$ defined on $\Gamma$ satisfies \eqref{ineq_fdecay} for some $0<\rho\le{1}/{2}$. Then, the associated maximal operator is bounded on $L^p$ for $p>{1}/{\rho}$.
\end{conjecture}
The conjecture was originally formulated by Stein in the case $\rho=1/2$, and later extended by Iosevich--Sawyer \cite{IS2} to $0<\rho\le 1/2$. Note that the conjecture only concerns $p>2$; when $p\le2$ the maximal estimates are influenced by various geometric features of $\Gamma$ beyond the Fourier decay (see \cite{NSW, IS2, ISS, HHY2018, BDIM2019,BIM}). 
For general hypersurfaces in the regime $p>2$, we refer works of Sogge--Stein \cite{SS1985} and the last author \cite{Oh2025}, which give partial results.
For a particular case $\rho=1/2$, however, the last author \cite{Oh2025} resolved the conjecture, which is Stein's conjecture.

The above conjecture asserts that \eqref{ineq_fdecay} is a sufficient condition for $L^p$ bounds of $M_\Gamma$. We can naturally ask that whether \eqref{ineq_fdecay} is also a necessary condition for such boundedness. Unfortunately, \eqref{ineq_fdecay} is far from being necessary, and in certain instances, it fails in strikingly drastic manner. 
 For analytic hypersurfaces in $\mathbb R^3$, Zimmermann \cite{Z} proved that $M_{\Gamma}$ is bounded on $L^p(\mathbb R^3)$ for $p>2$ whenever $\Gamma$ contains the origin, which is the only degenerate point of $\Gamma$.
 In such cases, therefore, maximal estimates still hold for $p>2$, even with extremely weaker Fourier decay.
This shows that, beyond Fourier decay, an additional geometric condition plays a role in the $L^p$ boundedness of maximal averages.
In fact,  the origin on $\Gamma$ in \cite{Z} is the simplest instance of what we will refer to as \emph{non-transversality}.

\begin{definition}\label{def_tran}
A point $x$ on $\Gamma$ is said to be \emph{transversal}, if the affine tangent space $x+T_x\Gamma$  to $\Gamma$ does not pass through the origin.
Otherwise, $x$ is called a \emph{non-transversal} point.

\end{definition}

In $\mathbb R^3$, Ikromov--Kempe--M\"uller \cite{IKM} considered hypersurfaces that consist only of transversal points. Under the transversality assumption, the authors proved Conjecture~\ref{conj_stein} for smooth, finite-type hypersurfaces in $\mathbb R^3$. 
Later, Greenblatt \cite{Gre2013, Gre2021} showed that the maximal estimates of \cite{IKM} remain valid even in the absence of transversality.
The maximal estimates of \cite{Gre2013, Gre2021}, despite accommodating non-transversal points, fall short of capturing the phenomena exhibited in Zimmermann's result \cite{Z}.

In this paper, we address this gap by formulating Conjecture \ref{conj_stein} in terms of explicit necessary and sufficient conditions when $\Gamma$ is an analytic hypersurface allowing non-transversal points. We also prove the reformulated conjecture for $\rho=1/2$. The key insight underlying this reformulation is that the relevant obstruction is intrinsically tied to non-transversality.

\subsection{Maximal estimates associated to non-transversality}

In $\bR^2$, Iosevich \cite{I} proved Conjecture~\ref{conj_stein} for curves of finite type. Later, Iosevich--Sawyer \cite{IS1} extended the result of Iosevich \cite{I} to the case of homogeneous hypersurfaces.
For the necessary part, it is proved in \cite{IS1} that $M_\Gamma$ is bounded on $L^p(\mathbb R^d)$ only if $d(\cdot , H)^{-1}\in L^{1/p}_\text{loc}(\Gamma,\mu)$ for all hyperplanes $H$ not passing through the origin, and thus hyperplanes tangent to $\Gamma$ at a non-transversal point are excluded. In \cite{IS2}, the authors conjectured that for $p>2$ the above necessary condition is in fact sufficient to ensure the $L^p$ boundedness of $M_\Gamma$.

The conjecture of Iosevich--Sawyer suggests that $M_\Gamma$ is bounded on $L^p$ for a wider range of $p$ when $\Gamma$ is localized near a non-transversal point. This can indeed be verified in certain specific cases, when $\Gamma$ is a curve in $\bR^2$ or a convex finite line type hypersurface in $\bR^d$, by the arguments in \cite{I,IS2} (see Appendix \ref{sec_scaling}). However,  for general hypersurfaces, no results are currently known that explain such improvements in the presence of non-transversal points.
Previously known results for non-transversal cases are limited to special classes of hypersurfaces (\cite{IS1,IS2}) or to the case when the non-transversal point is the origin (\cite{Z,LWZ}).
In this paper, we prove that $M_\Gamma$ is bounded on $L^p(\mathbb R^d)$ for every $p>2$,  whenever $\Gamma$ denotes a sufficiently small analytic hypersurface containing a non-transversal point.

Let $\gamma : [-2,2]^{d-1} \to \mathbb R$ be an analytic function and $\Gamma$ be the graph $\{(y, \gamma(y))  \in \mathbb R^d : y\in\mathbb [-2,2]^{d-1}\}$.
For a smooth function $\psi$ satisfying $\supp(\psi)\subset \{y\in\bR^{d-1}:|y|\le1/2\}$, 
let a measure $\mu[\psi]$ over $\Gamma$ be given by
\begin{align}
\label{def:mupsi}
	\langle \mu[\psi], f \rangle = \int_{\mathbb R^{d-1}} f (y, \gamma(y)) \psi(y)~\mathrm{d}y.
\end{align}
The dilation $\mu_t[\psi]$ by $t>0$ of $\mu[\psi]$ is defined as before. 
 Then we define
\[ M_\Gamma[\psi]f(x)\coloneq \sup_{t>0}|A_\Gamma[\psi]f(t,x)|\coloneqq \sup_{t>0}|\mu_t[\psi]\ast f(x)|. \]
The following theorem is our first main result on $M_\Gamma[\psi]$.
\begin{theorem}\label{thm_main}
    Let $d\ge2$ and $\Gamma$ be an analytic hypersurface. Suppose that a point $(y_{nt}, \gamma(y_{nt}))$ is non-transversal. If $\psi$ has a sufficiently small support around $y_{nt}$,
        then $M_{\Gamma}[\psi]$ is bounded on $L^p(\bR^d)$ for $p>2$.
\end{theorem}
Theorem~\ref{thm_main} gives a definitive result for the non-transversal case. In particular, it extends the result of Zimmermann \cite{Z} to every dimension and arbitrary non-transversal points.
 The range $p>2$ is sharp by considering cylindrical extension of a finite type plane curve. Moreover, together with  \cite[Corollary~1.8]{IKM}, Theorem~\ref{thm_main} resolves the conjecture of Iosevich--Sawyer \cite[Conjecture~1]{IS2} in the case of analytic surfaces in $\mathbb R^3$.  
The analytic condition in Theorem \ref{thm_main} is essential, as Zimmermann's example \cite{Z} shows that Theorem \ref{thm_main} cannot be extended to general smooth hypersurfaces. However, determining the precise class of smooth functions for which  Theorem \ref{thm_main} remains valid is still an open problem.

We note that the proof of Theorem~\ref{thm_main} relies neither on the resolution of singularity nor on the Newton diagram, which have been crucial tools in the literature \cite{IKM,Z, BDIM2019, Gre2021}. 
Instead, we decompose the surface $\Gamma$ into pieces according to the degeneracy  
of $\gamma$,
and by the $o$-minimal structure related to analytic functions, 
we exploit favorable structural properties within each decomposed piece.

In earlier works, the classical scaling argument has played an important role in verifying favorable properties of decomposed pieces.
The same idea also proved useful in \cite{IS2}, since a convex smooth function of finite line type could be approximated by a mixed homogeneous polynomial (see \cite{schulz1991}).
However, the scaling method becomes difficult to apply for general hypersurfaces.
In contrast, our method effectively replaces the role of scaling: it enables a direct analysis of each decomposed piece of an arbitrary analytic hypersurface, even in the absence of a scaling structure. We also note that when $\Gamma$ enjoys an inherent scaling structure, a simple proof of Theorem~\ref{thm_main} is available (see Appendix~\ref{sec_scaling}).

Another surprising consequence of the decomposition is a refined $L^2$ estimates for $M_\Gamma[\psi]$, which is generally unattainable without imposing the Fourier decay condition \eqref{ineq_fdecay} for $\rho=1/2$. A standard approach to proving the maximal estimates relies on regularity properties of $A_\Gamma[\psi]$ via the Sobolev embedding lemma. However, a direct application of this method fails to yield $L^2$ maximal estimates, since the $L^2$ regularity property of $A_\Gamma[\psi]$ is determined by \eqref{ineq_fdecay}. This difficulty is overcome by our decomposition according to the degeneracy of the hypersurface $\Gamma$. On each patch of $\Gamma$, the decomposition optimizes the regularity loss typically incurred through the Sobolev embedding lemma so that the loss and the localized $L^2$ regularity of the associated averaging operators are precisely balanced.
This compensation arises from non-transversality and the \L ojasiewicz inequality \cite{Loj1965} developed in real algebraic geometry. For a sufficiently large $p>2$, on the other hand, we adapt smoothing estimates of \cite{Oh2025} since each subsurface enjoys conical nature from the decomposition process.

\subsection{Reformulation of Iosevich--Sawyer--Stein conjecture}
Theorem~\ref{thm_main} reveals a striking dichotomy in the role of the geometry of hypersurfaces. Precisely, the local behavior of an analytic surface near a non-transversal point turns out to be irrelevant to the $L^p$ boundedness of the associated maximal operator whenever $p>2$. In contrast, near transversal points, the Fourier decay plays a crucial role in determining $L^p$ boundedness of maximal averages, as shown in previous works \cite{I,IS1, IS2, IKM}. In this sense, we re-formulate the Iosevich--Sawyer--Stein conjecture for analytic hypersurfaces.  
\begin{conjecture}\label{conj_iss_refined}
    Let $p_{cr}\ge2$ and $\Gamma$ be a compact analytic hypersurface. Then, the following statements are equivalent.
    \begin{enumerate}
        \item $M_{\Gamma}$ is bounded on $L^p(\bR^d)$ for all $p>p_{cr}$.
        \item For $\sigma\in C^\infty(\Gamma)$ and $p>p_{cr}$, there exists a constant $C>0$ such that
        $$
		|\cF[\sigma \, \mathrm{d}\mu](\xi)|\le C |\xi|^{-1/p},
	$$
    provided that $\sigma$ vanishes on an open set containing all non-transversal points of $\Gamma$.
    \end{enumerate}
\end{conjecture}

One may observe that in Conjecture~\ref{conj_iss_refined}, the Fourier decay is expected to be $1/p$ with $p>p_{cr}$, rather than just $1/p_{cr}$.
This is essential. By \cite{IM2011}, there is an analytic finite type surface in $\mathbb R^3$ whose measure $\mu$ satisfies a sharp decay estimate, 
\[
    |\widehat{\mathrm d\mu}(\xi)|\le C |\xi|^{-\frac{1}{2}} \log(2+|\xi|),\quad |\xi|\ge 1.
\]
Together with local smoothing estimates of \cite[Theorem 1.3]{Oh2025}, the associated maximal function is bounded on $L^p(\mathbb R^d)$ for all $p>2$.
From this example, $(2)$ of Conjecture~\ref{conj_iss_refined} is necessary for the refined formulation.

Thanks to Theorem \ref{thm_main}, Conjecture~\ref{conj_stein} and the implication $(2)\to(1)$ of Conjecture~\ref{conj_iss_refined} are equivalent for analytic hypersurfaces.
In particular, by the results of \cite{I,IS2,IKM}, the direction $(2)\to(1)$ of Conjecture~\ref{conj_iss_refined} holds for $\Gamma$ being a curve in $\bR^2$, a surface in $\bR^3$, or a convex hypersurface of finite line type in $\bR^d$.

The second main result of this paper contains a resolution of the relation $(1) \to (2)$ of Conjecture~\ref{conj_iss_refined}.
\begin{theorem}\label{thm_Fdecay}
    Let  $p_{cr}\ge1$ and $\Gamma$ be a compact analytic hypersurface.
    Suppose that $M_\Gamma$ is bounded on $L^p(\mathbb R^d)$ for $p>p_{cr}$ and $\sigma\in C^\infty(\Gamma)$ vanishes on an open set containing all non-transversal points of $\Gamma$.
    Then, there exists a constant $C_{p,\sigma}$ such that
    \[
        |\mathcal F  [\sigma \,\mathrm d\mu] (\xi)| \le C_{p,\sigma} |\xi|^{-1/p},\quad \forall p>p_{cr}.
    \]
\end{theorem}

To prove Theorem~\ref{thm_Fdecay}, we first return to the necessary condition of Iosevich--Sawyer \cite{IS2}:
\begin{align}\label{ineq_is_dist}
	\int_\Gamma d(x,H)^{-1/p} \chi(x)~\mathrm d\mu (x)<\infty,\quad p>p_{cr}
\end{align}
for all hyperplanes $H$ which do not pass through the origin. Note that $\chi$ denotes a cut-off function for the local integrability. 
We show that \eqref{ineq_is_dist} holds uniformly in the choices of $H$ 
under the assumption of $M_\Gamma$ being bounded on $L^p(\mathbb R^d)$. 
By a simple use of Chebyshev's inequality, the integrability of the distance functions yields a sublevel set estimate related to the defining function of $\Gamma$.
Then, we adapt a stationary set method of Basu--Guo--Zhang--Zorin-Kranich \cite{BGZZ} developed via model theory. The stationary set method helps us to control an oscillatory integral in terms of a \emph{mid-level} set estimate, which is a kind of sublevel set estimate. Thus, the overall argument proceeds through the implications (roughly) $(1)\to\eqref{ineq_is_dist}\to(2)$ of Conjecture~\ref{conj_iss_refined}.

As a corollary of Theorem~\ref{thm_Fdecay}, we prove Conjecture~\ref{conj_iss_refined} for $p_{cr}=2$.
\begin{corollary}
	Conjecture~\ref{conj_iss_refined} holds for $p_{cr}=2$.
\end{corollary}
Note that the $(1)\to (2)$ is obtained by Theorem~\ref{thm_Fdecay} with $p_{cr} =2$, and the $(2)\to(1)$ of Conjecture~\ref{conj_iss_refined} is obtained by \cite{Oh2025} and Theorem~\ref{thm_main}. For $(2)\to(1)$, one applies the local smoothing estimates \cite[Theorem 1.3]{Oh2025} to interpolate with trivial $L^2$ bounds, and obtains the range $p>2$.

\subsection*{Organization}
In Section~\ref{sec_pre}, we reduce Theorem~\ref{thm_main} to a refined $L^2$ estimate and smoothing estimates by means of a decomposition based on sizes of derivatives of $\gamma$. We also introduce $o$-minimal expansions of $\mathbb R$, which allow us to exploit useful properties of the decomposition.
Through Sections~\ref{sec_L2} and \ref{sec_smoothing}, we prove Theorem~\ref{thm_main} by verifying the refined $L^2$ estimates and smoothing estimates, respectively.
In Section~\ref{sec_fdecay}, we prove Theorem~\ref{thm_Fdecay} by revealing the relations among maximal estimates, sublevel set estimates, and Fourier decay.
In the appendix, we suggest a simple proof for Theorem~\ref{thm_main} in the case that hypersurfaces enjoy certain scaling structure.

\subsection*{Notation}
\begin{enumerate}
    \item We use both $\mathcal{F}$ and $^\wedge$ for the Fourier transform. 
    \item The notation $A\lesssim B$ means $A\leq CB$ for a harmless constant $C$.
    \item Let $\bB_r^N(p)$ be the ball of radius $r$ with center $p$ in $\bR^{N}$. We sometimes write $\bB_r^N$ for any ball of radius $r$ when the center is not determined.
    \item For a smooth function $\beta_0$ defined on $\bR$ satisfying $\beta_0 \equiv 1$ on $[-1,1]$ and $\beta_0 \equiv 0$ on $[-2,2]^\complement$, we set $\beta_1 = \beta_0(\cdot/2) - \beta_0(\cdot)$ so that $\sum_{\ell\in\bZ} \beta_1(s/2^\ell) =1$ for $s\neq 0$.
\end{enumerate}

\section{Preliminaries}\label{sec_pre}

Suppose $x_{nt}=(y_{nt},\gamma(y_{nt}))$ is a non-transversal point of $\Gamma$ and $\psi$ is supported on $\mathbb B^{d-1}_{1/2}(0)$. Without loss of generality, after a suitable rotation, we may assume that 
\begin{equation}
\label{gamma00}
   \gamma(y_{nt})=0, \quad \nabla\gamma(y_{nt})=0.
\end{equation}
 It is worth noting that when $x_{nt}$ is a transversal point, one may take $y_{nt}=0$ and assume that $\gamma(0)=c\neq0$, $\nabla\gamma(0)=0$ without loss of generality, however, this simplification is not possible at a non-transversal point.
Moreover, we assume that $\gamma$ is of finite type, in that, for each $y\in \bB_1^{d-1}(0)$, there exists a multi-index $\alpha=\alpha(y)$ with $|\alpha|\ge2$ such that $|\partial^\alpha\gamma(y)|\neq0$. Indeed, if $\gamma$ is not of finite type, then analyticity together with \eqref{gamma00} implies that $\gamma$ is identically zero. In this case, Theorem \ref{thm_main} follows immediately from the $L^p$-boundedness of the Hardy-Littlewood maximal function in $\bR^{d-1}$. 

By smoothness and finite type condition of $\gamma$, we have quantitative upper bounds on the derivatives of $\gamma$ and a lower bound on $k$-th derivatives of $\gamma$ for some $k\ge2$. More precisely, we may assume that there exists a constant $C_*$ such that
\begin{equation}
\label{upperbounds:gammma}
    \sup_{y\in\bB_1(0)}|\partial^\alpha\gamma(y)|\le C_*, \quad \text{ for all }|\alpha|\le N,
\end{equation}
for sufficiently large $N$.
After reducing $\supp(\psi)$ to a sufficiently small neighborhood of $y_{nt}$,
we may further assume that there exists an integer $k\ge2$ such that
\begin{equation}
\label{finitetype}
    \sum_{|\alpha|=k}|\partial^\alpha\gamma(y)|\ge C_*^{-1}, \quad \text{for }y\in \supp(\psi).
\end{equation}
When $k=2$, the $L^p$ boundedness of the corresponding maximal operator is well-known by Sogge \cite{Sogge}. Therefore, in the proof of the Theorem \ref{thm_main}, we assume that \eqref{gamma00}, \eqref{upperbounds:gammma} with sufficiently large $N$, and \eqref{finitetype} with $k\ge3$ hold.

In this section, we recall some useful results of real algebraic geometry and model theory, which allow us to decompose hypersurfaces according to the geometry around non-transversal points, and obtain practical properties of analytic functions.
As a consequence, Theorem~\ref{thm_main} is reduced to a study on localized operators under these settings.
\subsection{The \L ojasiewicz inequality and a decomposition of analytic hypersurfaces}
Here, we introduce our technical lemma in which the analyticity of a function is crucially used. 

\begin{lemma}[\cite{Loj1965}, \L ojasiewicz inequality]\label{lem:loj}
	Let $U$ be an open set in $\bR^{d-1}$ and $f : U \to \bR$ be an analytic function.
	For every point $x_0\in U$, there are a neighborhood $V$ of $x_0$, an exponent $\theta \in [1/2, 1)$, and a constant $C$ such that
	$$
		|f(x) - f(x_0)|^\theta \leq C |\nabla f (x)|,\quad \forall x \in V.
	$$
\end{lemma}

We take $x_0=y_{nt}\in \mathbb B_{1/2}^{d-1}(0)$ and $f(y)=\gamma(y)$. Considering \eqref{gamma00}, there exist a neighborhood $U_1$ of $y_{nt}$, an exponent $\theta_1\in[1/2,1)$ and a constant $C_1$ such that
\begin{align*}
	|\gamma(y)| \leq |\gamma(y) -\gamma(y_{nt})|^{\theta_1} \leq C_1 |\nabla \gamma(y)|
\end{align*}
holds on $U_1$. If we apply the \L ojasiewicz inequality on each component of $\nabla \gamma$, then there are a subset $U_2 \ni y_{nt}$ of $U_1$, an exponent $\theta_2\in[1/2, 1)$, and a constant $C_2$ such that for $y\in U_2$
\begin{align*}
	 |\nabla \gamma(y)| \leq |\nabla\gamma(y) -\nabla\gamma(y_{nt})|^{\theta_2} \leq C_2 |D^2 \gamma (y)|.
\end{align*}
Here, $|M|$ denotes the Hilbert-Schmidt norm for $M\in \mathbb R^{d\times d}$.
Therefore, it follows that for $y \in U_2$
\begin{align}\label{ineq_Loja}
	|\gamma(y)| \leq C_1 |\nabla \gamma(y)| \leq C_1 C_2 |D^2 \gamma (y)|.
\end{align}
Since we may choose the support of cut-off function $\psi$ sufficiently small, we assume that $\supp(\psi)\subset U_2$.
Thus, it suffices to prove Theorem \ref{thm_main}  under the assumption that \eqref{ineq_Loja} is valid on $\supp(\psi)$.

From now on, we give a decomposition of the cut-off function $\psi$ according to the size of derivatives of $\gamma$.
For this purpose, we define the following quantity:
\begin{align}
\nonumber
	\mathfrak{D}^m \gamma (y) := \Big( \sum_{|\alpha| = m} |\partial^\alpha \gamma(y)|^2 \Big)^{1/2}.
\end{align}
Note that $|\nabla \gamma(y)|$ and $|D^2\gamma(y)|$ are equivalent to
$$
	\Big(\sum_{|\alpha| = 1} |\partial^\alpha \gamma(y)|^2\Big)^{1/2},\quad 
	\Big(\sum_{|\alpha| = 2} |\partial^\alpha \gamma(y)|^2\Big)^{1/2},
$$
respectively. 
From \eqref{ineq_Loja}, the derivatives of $\gamma$ cannot vanish outside of the zeros of $\gamma$, 
and  $\mathfrak{D}^k\gamma(y) \sim1$ for all $y\in \supp(\psi)$ from the assumption \eqref{finitetype}.

Now we decompose $\supp(\psi)$ in terms of the size of $\mathfrak{D}^m\gamma$, $m=2, \dots, k$.
Let $\vec{\ell}=(\ell_2,\cdots,\ell_k) \in \bZ^{k-1}$. Then, we define $\psi_{\vec{\ell}}$ as following:
\begin{align}
\nonumber
	\psi_{\vec{\ell}\,}(y) = \psi(y) \prod_{m=2}^k \beta_1(\mathfrak{D}^m \gamma(y)/2^{\ell_m}).
\end{align}
By \eqref{upperbounds:gammma}, there exists a constant $C>0$ such that $\psi_\vl=0$ unless $\ell_m \leq C$ for all $m=2,\dots, k$. 
Thus, one has
\begin{equation}
    \label{decomp:psi}
	\psi(y) = \sum_{\vec{\ell}\,: \ell_m \leq C } \psi_{\vec{\ell}\,}(y),
\end{equation}
for $y\in\supp (\psi)$ satisfying $\mathfrak D^m\gamma(y)\neq0$ for all $m=2,\cdots,k$. One can check that the set $\cup_{m=2}^k\{y\in \supp(\psi):\mathfrak D^m\gamma(y)=0\}$ is of measure zero.

Since \eqref{decomp:psi} holds except a measure zero set, we have 
that $M_\Gamma[\psi] \leq \sum_{\vec{\ell}} M_\Gamma[\psi_{\vec{\ell}\,}]$. Hence it suffices to show the following result to prove Theorem \ref{thm_main} since $\ell_m$'s are bounded above by a fixed constant.
\begin{proposition}\label{prop_main}
For $p>2$, there exists a positive constant $c$ such that for all $\vl\in \mathbb Z^{k-1}$,
	$$\|M_\Gamma[\psi_{\vec{\ell}\,}] f\|_{L^p(\bR^d)} \lesssim 2^{c(\ell_2+\cdots+\ell_k)}\|f\|_{L^p(\bR^d)}.$$
\end{proposition}

Let $P_\lambda$ be a frequency projection operator onto $\{\xi\in\bR^d:\lambda\leq|\xi|\leq4\lambda\}$ given by $\widehat{P_\lambda g}(\xi) = \beta_1(|\xi|/\lambda) \widehat{g}(\xi)$.
To prove Proposition~\ref{prop_main}, it suffices to consider estimates for a local maximal function $\sup_{1<t<2} | A[\psi_{\vec{\ell}\,}] P_\lambda f |$. 
Then, Proposition~\ref{prop_main} is obtained via interpolation  and the Littlewood--Paley theory using the following two results. The argument is standard, so we just refer to \cite{Sch} or \cite[Lemma~3.4]{BHS}, for example. Sections~\ref{sec_L2} and \ref{sec_smoothing} are devoted to prove the following two propositions.

\begin{proposition}\label{prop_key}
	For $\lambda \geq 1$ and $\vl \in \mathbb Z^{k-1}$, we have
	$$
		\Big\|\sup_{1<t<2} |A_\Gamma[\psi_{\vec{\ell}\,}] P_\lambda f|\Big\|_{L^2(\bR^d)} \lesssim  \|f \|_{L^2(\bR^d)}.
	$$
\end{proposition}

\begin{proposition}\label{prop_smoothing}
	Let $\lambda\geq1$, $\vl \in \mathbb Z^{k-1}$, and $p\ge 4(k-1)$. 
	Then there exist positive constants $c$ and $\varepsilon$ such that
	$$
		\Big\|\sup_{1<t<2} |A_\Gamma[\psi_{\vec{\ell}\,}] P_\lambda f|\Big\|_{L^p(\bR^d)} \lesssim 2^{c(\ell_2+\cdots+\ell_k)}\lambda^{-\varepsilon}\|f \|_{L^p(\bR^d)}.
	$$
\end{proposition}
In the proof of Proposition \ref{prop_key} and Proposition \ref{prop_smoothing}, we may assume that 
\[
2^{\ell_k}\sim 1,
\]
since $\mathfrak D^k\gamma(y)\sim 1$ for all $y\in \supp(\psi)$. Note that $\psi_\vl=0$ whenever $\ell_k\le -C$ for some positive constant $C$. In this case the propositions become trivial.

In the rest of this section, we establish a size estimate for 
 $\supp(\psi_{\vec\ell})$, which plays a crucial role in the proofs of the two propositions above.
In \cite{Oh2025}, it was observed that the size of the set $\supp (\psi_{\vec\ell\,})\cap \mathbb B_{\epsilon}$ is bounded by
\begin{align}\label{def_delta}
\delta_{\vec\ell\,}:=\min_{m=2,\cdots,k-1}2^{\ell_m-\ell_{m+1}},
\end{align}
if $\epsilon$ is small enough.
For our purpose, we require an analogue of this size estimate that holds without any smallness assumption on $\epsilon$. For verifying this, we will make use of certain properties of $o$-minimal expansions of $\bR$, which we briefly review in the next subsection.

\subsection{$o$-minimal expansions of $\bR$}

In this subsection, we focus on two well-known expansions of the real ordered field $\bR=(\bR,(0,1),(+,\cdot),\le)$.
The first is
\[
\bR_{an}=(\bR,(0,1),(+,\cdot, \text{all  restricted   analytic   functions}),\le),
\]
which is known to be $o$-minimal. The second is
\[
\bR_{an,exp}=(\bR,(0,1),(+,\cdot, \text{all   restricted   analytic   functions},\exp),\le),
\]
whose $o$-minimality was established by \cite{DMM}.
Here, $\exp$ is the exponential symbol. A function $f:\bR^n\rightarrow \bR$ is called a restricted analytic function, if there exists an open set $U$ containing $[-1,1]^n$ and an analytic function $g:U\rightarrow \bR$ satisfying $f=g$ on $[-1,1]^n$ and $f=0$ otherwise.

We mainly consider functions which are \emph{definable} in the above two structures. We refer the readers to see \cite[Section 2]{BGZZ} for a precise definition of definability and $o$-minimal structures. In this paper, we only use the following basic properties: 
\begin{enumerate}
    \item All restricted analytic functions and polynomials are definable in $\bR_{an}$.
    \item A definable function $f$ in $\bR_{an}$ is definable in $\bR_{an,exp}$. 
    \item The collection of definable functions is closed under addition and multiplication. 
    \item If $f$ is definable, then so is $g=\mathds1_{\{f\ge0\}}$, where $\mathds1_U$ denotes the characteristic function supported on $U$.
\end{enumerate}
In addition, by \cite{CM}, we have the following proposition.
\begin{proposition}[\cite{LO}, Proposition 3.6]\label{prop:LO}
    Let $f:\bR^n\times\bR^m\rightarrow \bR$ be definable in $\bR_{an}$. If $f(\cdot,y)$ is integrable for all $y$, then
    \[
    y\in \bR^m\mapsto g(y)=\int_{\bR^n} f(x,y)~\mathrm dx 
    \]
    is definable in $\bR_{an,exp}$.
\end{proposition}

By these properties, we can give examples of definable functions. For instance, $\partial^\alpha\gamma \mathds 1_{[-1,1]^{d-1}}$ is definable in $\bR_{an}$ for every multi-index $\alpha \in\bN_0^{d-1}$. Additionally,
\[
(w,r,\delta)\mapsto \int_{\mathbb R^{d-1}}\mathds1_W(y,w,r,\delta)~\mathrm dy
\]
is definable in $\bR_{an,exp}$ where 
\[
W=\{(y,w,r,\delta)\in [-1,1]^{d-1}\times \bR^{d-1}\times \bR\times\bR:|\gamma(y)+w\cdot y-r|\le\delta\}.
\]

Such functions, which are definable in $o$-minimal structure, have the important property that the number of monotonicity changes is uniformly bounded. More precisely, the following holds.
\begin{proposition}[\cite{BGZZ}, Proposition 2.8]\label{prop:BGZZ}
    Let $h:\bR\times \bR^n\rightarrow \bR$ be a definable function in $\bR_{an,exp}$ and $N(y)$ be the number of times that a map $r\mapsto h(r,y)$ changes monotonicity. Then $\sup_yN(y)<\infty$ holds.
\end{proposition}

This proposition indicates that definable functions share certain properties with polynomials. In \cite{Oh2025}, a restriction on a small ball is needed to estimate $|\supp (\psi_{\vec\ell \,})|$, since $\gamma$ can then be approximated by a polynomial. In contrast, the favorable properties of definable functions allow us to estimate the size of $\supp (\psi_{\vec\ell \,})$ without any restriction on the support, as stated in the following proposition. This proposition directly implies that
 \begin{equation}
 \label{ineq:suppsize}
    |\supp (\psi_{\vec\ell\,})|\lesssim \delta_{\vec\ell\,}    
 \end{equation}
holds, which will be frequently used later. 
\begin{proposition}\label{prop_supp_size}
    Let $\gamma$ be an analytic function on $[-2,2]^{d-1}$ and $m\in \bN$. For all positive constants $c_m$, we have
   \[
   |\{y\in[-1,1]^{d-1}:\mathfrak D^m\gamma(y)\le c_m,\mathfrak D^{m+1}\gamma(y)\ge c_{m+1}\}|\lesssim c_{m}/c_{m+1}.
   \]
   Here, an implicit constant depends only on $m$ and $\gamma$.
\end{proposition}  
\begin{proof}
Define a set
\[
U:=\{y\in[-1,1]^{d-1}:\mathfrak D^m\gamma(y)\le c_m,\mathfrak D^{m+1}\gamma(y)\ge c_{m+1}\}.
\]
    Since $c_{m+1}^{-1}\mathfrak D^{m+1}\gamma\ge 1$ on $U$, one can check that
    \[
    |U\,|\le \int_{U} c_{m+1}^{-1}\mathfrak D^{m+1}\gamma(y)~\mathrm{d}y.
    \]
    Considering $\mathfrak D^{m+1}\gamma\sim \sum_{|\alpha|=m+1}|\partial^\alpha\gamma| $, it is enough to show that \[
        \int_{U}|\partial^\alpha\gamma|\lesssim c_m,\quad\text{for}\quad|\alpha|=m+1.
    \]
By the fundamental theorem of calculus and Fubini's theorem, the above inequality follows if we prove that there exists a number $N>1$ depending only on $\gamma$ such that for all $j=1,\cdots,d-1$, $|\alpha|=m+1$, and $y\in[-1,1]^{d-1}$, the function  
\begin{equation}
\label{function:jalpha}
    r\in [-1,1] \mapsto  \partial^\alpha\gamma(y_1,\cdots,y_{j-1},r,y_{j+1},\cdots,y_{d-1})
\end{equation}
changes its sign at most $N$ times.

To show this, we use Proposition \ref{prop:BGZZ}. For fixed $j$ and $\alpha$ satisfying $|\alpha|=m$, let $h:\bR^{d-1}\rightarrow \bR$ be a function satisfying $h=\partial^\alpha\gamma$ on $[-1,1]^{d-1}$ and $h=0$ otherwise. Define $\tilde h:\bR\times\bR^{d-2}\rightarrow\bR$ by 
\[
\tilde h(r,z)=h(z_1\cdots,z_{j-1},r,z_{j},\cdots,z_{d-2}).
\]Then $\tilde h$ is a restricted analytic function, so it is a definable function in $\bR_{an}$. 
By Proposition \ref{prop:BGZZ}, we can conclude that the number of changes in monotonicity of the function $r\mapsto \tilde h(r,z)$ is bounded by a constant depending only on $\gamma$. Thus, there exists $N_{j,\alpha}>1$ such that for all $y\in[-1,1]^{d-1}$, the function \eqref{function:jalpha} changes its sign at most $N_{j,\alpha}$ times. Since the number of $j$ and $\alpha$ is finite, we get the desired estimate.
\end{proof}

\section{Proof of Proposition~\ref{prop_key}: $L^2$ bounds}\label{sec_L2}

By the Fourier inversion, we rewrite $A_\Gamma [\psi_{\vec{\ell}}\,]P_\lambda f$ as
\begin{align}
\nonumber
	A_\Gamma [\psi_{\vec{\ell}\,}]P_\lambda f 
	= &(2\pi)^{-d} \int_{\bR^d} m_{\vec{\ell}\,}(t\xi) \widehat{P_\lambda f}(\xi)  \mathrm{e}^{i x\cdot\xi}~\mathrm{d}\xi, \\
 \nonumber
   m_{\vec{\ell}\,}(t\xi) =& \int \mathrm{e}^{-i t \langle\Gamma(y), \xi \rangle}\psi_{\vec{\ell}\,}(y) \beta_1(|\xi|/\lambda)~\mathrm{d}y.
\end{align}
To show Proposition \ref{prop_key}, we utilize the Sobolev embedding lemma together with suitable multiplier bounds. Heuristically, replacing a supremum by an $L^2$-norm through the Sobolev embedding yields a $\lambda^{1/2}$-loss. In our situation, however, the loss would be refined since $\gamma$ is bounded above by $2^{\ell_2}$ in $\supp (\psi_{\vec\ell\,})$, which is not available in general. Indeed, one can observe that 
\begin{equation}
\label{partialtmell}
    \partial_t(m_{\vec{\ell}\,}(t\xi))=-i\langle\Gamma(y), \xi \rangle m_{\vec{\ell}\,}(t\xi),
\end{equation}
and the term $\langle\Gamma(y), \xi \rangle$ is typically much smaller than $\lambda$ when $\xi$ lies near the normal directions to $\Gamma$. This observation is crucial, as it leads to a significant improvement in the analysis. More precisely, the proof of Proposition \ref{prop_key} is reduced to $L^\infty$ bounds of both $m_{\vec{\ell}\,}(t\xi)$ and $\partial_t (m_{\vec{\ell}\,}(t\xi))$ by using the following Sobolev embedding lemma, which is well-known in the literature.
\begin{lemma}\label{sobo}
    Let $F\in C^1(\mathbb R)$ and $p>1$. Then
    \[
    \sup_{1<t<2}|F(t)|^p \le |F(1)|^p + p \Big(\int_1^2 |F(t)|^p\,\mathrm dt\Big)^{\frac{p-1}p} \Big(\int_1^2 |\partial_t F(t)|^p\,\mathrm dt\Big)^{1/p}.
    \]
\end{lemma}
We also refer to \cite[p.499, Lemma~1]{Stein} for a different version.
\begin{proof}
By the fundamental theorem of calculus, one has
\[
    |F(t)|^p \le |F(1)|^p + p \int_1^2 |F(t)|^{p-1}|\partial_t F(t)|\,\mathrm dt.
\]
The desired inequality follows by H\"older's inequality.
\end{proof}

Lemma \ref{sobo} shows that
\[
	\Big\|\sup_{1<t<2} |A_\Gamma[\psi_{\vec{\ell}\,}] P_\lambda f|\Big\|_{L^2(\bR^d)}
    \lesssim 
    \| A_\Gamma[\psi_{\vec{\ell}\,}] P_\lambda f (1, \cdot)\|_{L^2_{x}}
    +\prod_{i=0,1}\|\zeta\partial_t^iA_\Gamma[\psi_{\vec{\ell}\,}] P_\lambda f\|_{L^2_{t,x}}^{1/2},
\]
where $\zeta $ denotes a smooth cutoff function supported in $[1/2, 4]$ defined by
\[
\zeta(t)=\beta_1(t)+\beta_1(2t).
\]
Note that the first term in the RHS is trivially bounded by $\|f\|_2$. 
Therefore, by  Plancherel's theorem, 
Proposition \ref{prop_key} would follow if one could show 
\begin{equation}
\nonumber
    \Vert m_{\vec{\ell}\,}(t\cdot)\beta_1(|\cdot|/\lambda)\Vert_{L^\infty}\cdot\Vert \partial_t m_{\vec{\ell}\,}(t\cdot )\beta_1(|\cdot|/\lambda)\Vert_{L^\infty}\lesssim 1.
\end{equation}
However, this estimate does not hold in general, since the factor $\langle\Gamma(y), \xi \rangle$ in \eqref{partialtmell} may become large when $\xi$ is far from the normal direction to $\Gamma$. This indicates that the $L^\infty$-norm of $m_{\vec{\ell}\,}(t\cdot)$ and $\partial_tm_{\vec{\ell}\,}(t\cdot)$ may be attained by different $\xi$, considering the fact that the $L^\infty$-norm of $m_{\vec{\ell}\,}(t\cdot)$ is typically achieved when $\xi$ lies in a small neighborhood of the normal vector to $\Gamma$. Thus, it is natural to decompose $\xi$ in a suitable way to obtain the desired multiplier bounds. To this end, we decompose $\xi$ suitably into two pieces and establish the corresponding bounds for each piece, which gives an analogue of the above estimate. The decomposition is fulfilled according to the angle from the line spanned by $(0,\cdots,0,1)$, which is the normal vector to $\Gamma$ at $y=y_{nt}$.

Let $\xi=(\xi',\xi_d)\in \bR^{d-1}\times \bR$ and $\partial_i=\partial_{y_i}$.
Suppose that $\lambda\le |\xi|\le 4\lambda$ and $|\xi'|\geq 20C_2\lambda2^{\ell_2}$ where $C_2$ is the constant in \eqref{ineq_Loja}, then there exists an $1\leq i\leq d-1$ such that 
\begin{align}\label{lbound2}
	|\partial_{i} (\langle \Gamma(y), \xi\rangle)| = |\xi_i + \partial_i\gamma(y)\xi_d|>4C_2\lambda 2^{\ell_2}/\sqrt{d-1}.
\end{align}
Indeed, we have
\[|\xi'+\xi_d\nabla\gamma(y)|\geq |\xi'|-|\xi_d\nabla\gamma(y)|\geq 4C_2\lambda2^{\ell_2},\]
so that we can find a component satisfying \eqref{lbound2}. 
Therefore, in the case $|\xi'|\geq 20C_2\lambda2^{\ell_2}$, one can apply integration by parts together with the lower bound \eqref{lbound2} to improve multiplier bounds. 
From the observation, we split the frequency variable $\xi$ into two parts, by using
$P_c$ and  $P_o$ given by
\begin{align*}
   &\widehat{P_c g}(\xi) = \beta_0\Big(\frac{|\xi'|}{20C_2\lambda 2^{\ell_2}}\Big) \widehat{g}(\xi),\\
   &\widehat{P_o g}(\xi) = \Big(1 - \beta_0\Big(\frac{|\xi'|}{20C_2\lambda 2^{\ell_2}}\Big) \Big) \widehat{g}(\xi).
\end{align*}

Now, decompose $ A_\Gamma[\psi_{\vec{\ell}\,}]P_\lambda$ into two parts
\[
	 A_\Gamma[\psi_{\vec{\ell}\,}]P_\lambda =  A_\Gamma[\psi_{\vec{\ell}\,}]P_\lambda P_c +  A_\Gamma[\psi_{\vec{\ell}\,}]P_\lambda P_o.
\]
By the triangle inequality, to prove Proposition \ref{prop_key}, it is enough to show that 
\begin{equation}
\label{ineq:pc}
\Big\|\sup_{1<t<2} |A_\Gamma[\psi_{\vec{\ell}\,}] P_\lambda P_c f|\Big\|_{L^2(\bR^d)}\lesssim \|f\|_{L^2(\mathbb R^d)}
\end{equation}
and 
\begin{equation}
    \label{ineq:po}
    \Big\|\sup_{1<t<2} |A_\Gamma[\psi_{\vec{\ell}\,}] P_\lambda P_o f|\Big\|_{L^2(\bR^d)}\lesssim \|f\|_{L^2(\mathbb R^d)}.
\end{equation}

\subsection{Proof of \eqref{ineq:pc}}
We first consider \eqref{ineq:pc}. In this case we have $|\xi'| < 40C_2\lambda2^{\ell_2}$,  hence for $y\in \supp (\psi_{\vec\ell})$, we have that 
\[
|\langle\Gamma(y), \xi \rangle|\lesssim \lambda2^{\ell_2},
\]
due to the \L ojasiewicz inequality \eqref{ineq_Loja}.
Considering Lemma \ref{sobo} and \eqref{partialtmell}, the inequality \eqref{ineq:pc}
follows directly from the lemma below.
\begin{lemma}\label{lem_l2bound1}
Let $1/2\le t\le 4$. Then the following inequality holds when $\xi$ satisfies $|\xi|\sim\lambda$ and $|\xi'|\le 40C_2\lambda 2^{\ell_2}$; 
\[
|m_{\vec\ell\,}(t\xi)  |\lesssim (\lambda 2^{\ell_2})^{-1/2}.
\]
\end{lemma}

\begin{proof}
We first observe properties of $\gamma$ and $ \psi_{\vec\ell\,}$. In $\supp(\psi_{\vec\ell\,})$, we particularly have
\[
	2^{\ell_2}\leq \mathfrak{D}^2\gamma(y) \leq 4\cdot 2^{\ell_2}.
\]
Thus, for each $y$, there is a direction $v = v(y)$ such that $D_v^2 \gamma(y) \ge d^{-1}2^{\ell_2}$, where $D_v$ denotes the directional derivative along $v$.   On the other hand, since $\mathfrak D^m\gamma$ is smooth on a small neighborhood of $\supp(\psi_\vl)$, we can choose a small number $c>0$ such that for $2\le m \le k$,
\begin{equation}
\label{ineq:d3gamma}
    \mathfrak D^m\gamma(y)\le 8\cdot2^{\ell_m}, 
\end{equation}
if $\dist (y,\supp(\psi_\vl))\le c\delta_{\vl} $. Moreover, we can choose such a constant $c$ independent of $\vl$. Indeed, let $|y-z|\le c\delta_\vl$ for some $z\in \supp(\psi_\vl)$. Then, from \eqref{upperbounds:gammma}, $2^{\ell_k} \sim1$, and the fact that 
\[
|\partial_j\mathfrak D^m\gamma(y)|\le \Big(\sum_{|\alpha|=m}|\partial_j\partial^\alpha\gamma(y)|^2\Big)^{1/2}\le \mathfrak D^{m+1}\gamma(y),
\]
the Mean value theorem  gives that $\mathfrak D^{k}\gamma(y)\le 4\cdot 2^{\ell_k}+c\delta_{\vl}C_*k^d \le(4+c(C_* k^d)^4)2^{\ell_{k}}$. Here, we use $2^{\ell_{k-1}} \le C_* k^d$ and $(C_*k^d)^{-1}\le 2^{\ell_k} \le C_* k^d$. By the same reason and using  $2^{\ell_k}\delta_\vl\le 2^{\ell_{k-1}}$, the estimate on $\mathfrak D^{k}\gamma$ implies that $\mathfrak D^{k-1}\gamma(y)\le (4+c(4+c (C_* k^d)^4))2^{\ell_{k-1}}$. Thus, inductively, it follows that for $2\le m\le k$,
\[
\mathfrak D^m\gamma(y)\le \Big(c^{k+1-m} (C_* k^d)^4+4\sum_{j=0}^{k-m}c^{j}\Big)2^{\ell_m},
\]
which allows us to choose $c$ depending only on $C_*,k,d$ satisfying \eqref{ineq:d3gamma}.

Now, we decompose the support of $\psi_{\vec{\ell}\,}$ into balls of radius $c \delta_{\vec{\ell}\,}$ centered at $z\in c\delta_{\vec{\ell}\,} \mathbb{Z}^{d-1}$. If $c$ is sufficiently small, \eqref{ineq:d3gamma} with $m=3$ implies that for each $z$ satisfying $\supp(\psi_\vl)\cap\mathbb B_{c\delta_\vl}(z)\neq \emptyset$, there exists $z'\in \supp(\psi_\vl)\cap \mathbb B_{c\delta_\vl}(z)$ such that
\begin{align}\label{lbound_z}
	|D_{v(z')}^2 \gamma(y)| > \frac{2^{\ell_2}}{4d}\quad \forall y \in \mathbb{B}_{c \delta_{\vec{\ell}}}(z).
\end{align}
Similarly,
\eqref{ineq:d3gamma} yields for $y\in \mathbb B_{c\delta_\vl}(z)$,
\begin{align*}
    |\mathfrak{D}^m\gamma(y)| \sim 2^{\ell_m},\quad m=2, \dots, k-1.
\end{align*}
Note that we have $\mathfrak D^k \gamma (y) \sim 2^{\ell_k}$ for all $y \in \supp(\psi)$.

Thus, by Proposition~\ref{prop_supp_size}, the cardinality of such balls of radius $c\delta_{\vec{\ell}\,}$ is bounded by  constant times $\delta_{\vec{\ell}\,}^{2-d}$.
That is, 
\begin{align}\label{250401_2330}
	\# \Big( \Big\{\bB_{c\delta_{\vl}}(z) : z\in c \delta_{\vec{\ell}\,} \mathbb{Z}^{d-1},\,\, \dist(z, \supp(\psi_{\vec{\ell}\,})) \leq c\delta_{\vec{\ell}\,}\Big\} \Big) \lesssim \delta_{\vec{\ell}\,}^{2-d}.
\end{align}

Let $\{\varphi_z\}_z$ be a partition of unity subordinated to the family of balls $\{\mathbb{B}_{c \delta_{\vec{\ell}\,}}(z)\}$.
Then for each $z$, by \eqref{lbound_z}, scaling, and the van der Corput lemma along $v(z')$ for some fixed $z' \in \mathbb{B}_{c \delta_{\vec{\ell}\,}}(z)$, one has
\begin{align*}
	\Big| \int_{\bR^{d-1}} \mathrm{e}^{i t \Gamma(y)\cdot \xi} \psi_{\vec{\ell}\,}(y) \varphi_z(y)~\mathrm{d}y \Big|
	\lesssim (\lambda 2^{\ell_2} \delta_{\vec{\ell}\,}^2)^{-1/2} \delta_{\vec{\ell}\,}^{d-1}.
\end{align*}
Indeed, we first perform the change of variables $y\to \delta_{\vec{\ell}\,}\,y$ so that $|D_{v(z')}^2(\gamma(\delta_{\vec{\ell}\,}y))| \sim 2^{\ell_2} \delta_{\vec{\ell}\,}^2$. Hence, we have $|D_{v(z')}^2 t\Gamma(\delta_{\vec{\ell}\,} \,y)\cdot\xi|\gtrsim \lambda 2^{\ell_2} \delta_{\vec{\ell}\,}^2$ and $$|\partial_i (\psi_\vl\,(\delta_\vl \, y)\varphi_z(\delta_\vl\, y))| \lesssim \sum_{m=2}^k \mathfrak D^{m+1}\gamma(y) 2^{-\ell_m} \delta_\vl \lesssim 1. $$
The desired inequality follows by the van der Corput lemma and the scaling constant $\delta_{\vec{\ell}\,}^{d-1}$.
Together with the cardinality \eqref{250401_2330}, it follows that
\[
	\Big| \int_{\bR^{d-1}} \mathrm{e}^{i t \Gamma(y)\cdot \xi} \psi_{\vec{\ell}\,}(y)~\mathrm{d}y \Big|
	\lesssim  (\lambda 2^{\ell_2})^{-1/2},
\]
which is the desired bound on the multiplier.
\end{proof}

\subsection{Proof of \eqref{ineq:po}}
Now, we focus on the estimate \eqref{ineq:po}. In this case, the size of $|\langle\Gamma(y), \xi \rangle|$ depends on the size of $|\xi'|$.
To handle this effectively,  we further decompose $P_o$ with respect to the size of $\xi'$. 
Define the operator $P_{o,j}$ by
\[
\widehat{ P_{o, j}g}(\xi) = \beta_1\Big(\frac{|\xi'|}{20C_2 \lambda 2^{j}}\Big) \widehat{g}(\xi).
\]
Then by the following identity,
\begin{align}
\nonumber
    \beta_1(|\xi|\lambda^{-1}) \Big(1 - \beta_0\Big(\frac{|\xi'|}{20C_2\lambda 2^{\ell_2}}\Big) \Big) =\beta_1(|\xi|\lambda^{-1})  \sum_{j \ge \ell_2} \beta_1\Big(\frac{|\xi'|}{20C_2 \lambda 2^{j}}\Big),
\end{align}
one easily checks that $P_o=\sum_{\ell_2\le j\le 0}P_{o,j}$.
Note that $\ell_2$ is a negative number and the cardinality of $j$ is $|\ell_2|$.
By this decomposition and the triangle inequality, \eqref{ineq:po} follows once we show that there exists a constant $c>0$ such that 
\begin{equation}
\label{ineq:poj}
\|\sup_{1\le t\le 2} |  A_\Gamma[\psi_{\vec{\ell}\,}]P_\lambda P_{o,j} f |\|_{L^2(\mathbb R^d)}\lesssim 2^{c\ell_2} \|f\|_{L^2(\mathbb R^d)} \text{ for all }\ell_2\le j\le0.    
\end{equation}

When $1<|\xi'|/20C_2\lambda 2^j<4$ for an $\ell_2 \le j\leq 0$,  since $|\xi_d\gamma(y)|\lesssim 
\lambda 2^{\ell_2}$, one can see that 
\begin{equation}
\nonumber
    |\langle \Gamma(y),\xi\rangle|\lesssim \lambda 2^j \quad  \text{ for }y\in\supp(\psi_{\vec\ell\,}).
\end{equation} 
As in the case of $A[\psi_{\vec{\ell}\,}]P_{\lambda}P_c$, the desired estimate \eqref{ineq:poj} is a direct consequence of the following lemma.

\begin{lemma}
Let $\ell_2\le j\le0$ and $1/2\le t\le 4$. Then there exists a constant $c>0$ such that 
	\[
    |m_{\vec\ell\,}(t\xi)|\lesssim 2^{c\ell_2}(\lambda 2^j)^{-1/2},
    \]
    when $\lambda \le|\xi|\le 4\lambda$ and $20C_2\lambda2^j\le|\xi'| \le  80C_2\lambda 2^j$.
\end{lemma}

\begin{proof}
We apply integration by parts, using a lower bound on the gradient of the phase $\xi'\cdot y+\xi_d\gamma(y)$. First, we check the lower bound on the phase.
	Since we have
    \[
20C_2\lambda2^j<|\xi'| < 80C_2\lambda 2^j,\quad \text{and}\quad |\partial_i\gamma(y)\xi_d|\leq 16C_2\lambda2^{\ell_2} \le 16C_2 \lambda 2^j,
    \]
     it holds that $| \xi_1 +\partial_1 \gamma (y)\xi_d |\sim \lambda 2^j$ after a suitable rotation. 
	Applying integration by parts and the support condition $|\supp (\psi_{\vec\ell\,})|\lesssim \delta_{\vec\ell\,}$, one has
	\begin{align}
    \nonumber
		|m_{\vec{\ell}}\,(t\xi)| 
        \lesssim  \min\{ (\lambda 2^j)^{-1}, \delta_{\vec{\ell}\,}\}\le (\lambda 2^j)^{-1/2}\delta_{\vec\ell\,}^{1/2}.
	\end{align}
   The integration by parts produces, in addition to the main term $(\lambda 2^j)^{-1}$, an extra factor $\delta_{\vec{\ell}\,}^{-1}$ due to a derivative falling on $\psi_{\vec{\ell}}$. The factor $\delta_{\vec{\ell}\,}^{-1}$, however, is compensated by the support size $|\supp(\psi_{\vec\ell\,})| \lesssim \delta_{\vec{\ell}\,}$.
    Thus, the desired inequality holds by the observation that $\delta_{\vec\ell\,}^{k-2}\le 2^{\ell_2-\ell_3}\times\cdots\times2^{\ell_{k-1}-\ell_k}\sim 2^{\ell_2}$. 
\end{proof}

\section{Proof of Proposition~\ref{prop_smoothing}: Smoothing estimates}\label{sec_smoothing}

We prove Proposition~\ref{prop_smoothing} by considering two cases: 
\begin{itemize}
\item $2^{\ell_m}<\lambda^{-1/2}$ for some $m=2, \dots, k-1$.
\item $2^{\ell_m}\geq \lambda^{-1/2}$ for all $m$. 
\end{itemize}
In the first case, $\supp(\psi_{\vec{\ell}}\,)$ is sufficiently small, so we can get sufficient decay in the $L^p$ estimates for the frequency localized maximal function, in terms of the size of $\supp(\psi_{\vec{\ell}}\,)$ for any $p\geq2$.
On the other hand, it is required to perform a delicate frequency decomposition in the second situation, which eventually leads us to desired smoothing estimates.

\subsection{$2^{\ell_m}<\lambda^{-\frac12}$ for some $m=2,\dots, k-1$}

Recall that \eqref{def_delta} and \eqref{ineq:suppsize}.
Suppose $2^{\ell_{m_0}}<\lambda^{-\frac{1}{2}}$ for some $2\leq m_0\leq k-1$. Then, we have
\[ \delta_{\vec{\ell}\,}\leq (2^{\ell_{m_0} - \ell_{m_0+1}}\times\cdots\times 2^{\ell_{k-1} - \ell_k})^{\frac{1}{k-m_0}}. \]
Since $2^{\ell_k} \sim 1$, we have 
\begin{align}\label{sscale_size}
    | \supp(\psi_{\vec{\ell}\,})| \lesssim \lambda^{-\frac{1}{2(k-m_0)}}\leq \lambda^{-\frac{1}{2(k-2)}}.
\end{align}
Using $2^{\ell_k} \sim 1$ again, we also have
\begin{align}\label{sscale_sumf}
\begin{split}
	\delta_{\vec{\ell}\,} 
	\leq 2^{c (\ell_2 - \ell_3)} 2^{2c (\ell_3 - \ell_4)} \times \cdots\times 2^{(k-2)c(\ell_{k-1} - \ell_k)}
	\lesssim   2^{c(\ell_2 + \cdots + \ell_{k})} ,
\end{split}
\end{align}
where $c+2c+\cdots+(k-2)c=1$.
Thus, we have $c = \frac{2}{(k-1)(k-2)}$. On the other hand, by Minkowski's inequality, we obtain
\begin{align*}
    \| \zeta A_\Gamma[\psi_{\vec{\ell}\,}]P_\lambda f\|_{L^p(\mathbb{R}^{d+1})}
    \lesssim  &\,| \supp(\psi_{\vec{\ell}\,})| \, \| f\|_{L^p(\mathbb{R}^{d})},
\end{align*}
 where $\zeta(t)=\beta_1(t)+\beta_1(2t)$ as before.
Averaging \eqref{sscale_size} and \eqref{sscale_sumf}, we get
\begin{align}
\nonumber
    \| \zeta A_\Gamma[\psi_{\vec{\ell}\,}]P_\lambda f\|_{L^p(\mathbb{R}^{d+1})}
    \lesssim \lambda^{-\frac{1}{4(k-2)}} 2^{ \frac{c}{2}(\ell_2 + \cdots + \ell_{k})} \|f\|_{L^p(\mathbb{R}^{d})},
\end{align}
for any $p\geq1$. Note that Lemma~\ref{sobo} essentially incurs $\lambda^{1/p}$-loss when one handles the $L^p$ norm of local maximal average 
\[\sup_{1<t<2} |A_\Gamma[\psi_{\vec{\ell}\,}]P_\lambda f |.\] 
Hence, by the smoothing estimate and Lemma \ref{sobo}, we obtain
\begin{align}
\nonumber
    \Big\| \sup_{1<t<2} |A_\Gamma[\psi_{\vec{\ell}\,}]P_\lambda f | \Big\|_{L^p(\mathbb{R}^d)} \lesssim \lambda^{-\varepsilon(p)} 2^{\frac c2 (\ell_2 +\cdots + \ell_k)} \|f\|_{L^p(\mathbb{R}^{d})},
\end{align}
for some $\varepsilon(p)>0$ when $p>4(k-2)$, which gives Proposition \ref{prop_smoothing} for this case.

\subsection{$2^{\ell_m} \geq \lambda^{-\frac12}$ for all $m=2, \dots, k-1$}

By the Fourier transform in the $t$-variable, we write $A_\Gamma[\psi_{\vec{\ell}\,}]f(x, t)$ as
\[
	(2\pi)^{-d-1}\iint_{\mathbb{R}^{d+1}} \Big( \iint_{\mathbb{R}^d} \mathrm{e}^{-i s(\langle \Gamma(y), \xi\rangle + \tau)}\psi_{\vec{\ell}\,}(y) ~\mathrm{d}y\mathrm{d}s \Big) \widehat{f}(\xi) \mathrm{e}^{i(\langle x, \xi\rangle + t\tau)}~\mathrm{d}\xi \mathrm{d}\tau.
\]
By abusing notation, for a general symbol $\mathfrak{a} = \mathfrak{a}(y,s,\xi, \tau) \in C_c^\infty(\mathbb{B}^{d-1}(0,1) \times (1/2,4) \times \bR^{d}\times \bR)$, we denote
\[
	A_\Gamma[\mathfrak{a}]f(x,t)
	=\iint_{\mathbb{R}^{d+1}} \Big( \iint_{\mathbb{R}^d} \mathrm{e}^{-i s(\langle \Gamma(y), \xi\rangle + \tau)}\mathfrak{a}~\mathrm{d}y\mathrm{d}s\Big) \widehat{f}(\xi) \mathrm{e}^{i(\langle x, \xi\rangle + t\tau)}~\mathrm{d}\xi \mathrm{d}\tau.
\]
Then, for $\mathfrak a_\vl(y,s,\xi,\tau)=\psi_\vl\,(y)\zeta(s)$, we have $\zeta A_\Gamma[\psi_{\vec{\ell}\,}]f=A_\Gamma[\mathfrak a_\vl]f$.

Recall that we are considering $ A_\Gamma[\psi_{\vec\ell}]P_\lambda f(x,t)$ so that $|\xi|\sim \lambda$. Thus, if $\tau>C_\circ\lambda$ for some constant $C_\circ>0$, then we have $|\langle\Gamma(y),\xi\rangle +\tau|>\lambda$. Hence, the corresponding kernel satisfies a sufficiently large decay, which allows arbitrary smoothing estimates for $p\geq1$. 
Thus, we restrict $\tau$ into $[-C_\circ\lambda, C_\circ\lambda]$, and it suffices to consider  $A_\Gamma[\mathfrak{a}_{\vec{\ell}, \lambda}]f$, where
\begin{equation}
    \label{def:allambda}
    \mathfrak{a}_{\vec{\ell}, \lambda}\, (y, s, \xi, \tau) =  \psi_{\vec{\ell}\,}(y)  \zeta(s) \beta_1(\lambda^{-1}|\xi|) \beta_0((C_{\circ}\lambda)^{-1}\tau).
\end{equation}
By Lemma \ref{sobo}, Proposition~\ref{prop_smoothing} is reduced to the following result. 
\begin{proposition}\label{prop_smoothing_reduc}
	For $p>4(k-1)$ and $\varepsilon>0$, there are constants $C>0$ and $c >0$, independent of $\lambda$ and $\vl$, such that
	\[
		\| A_\Gamma[\mathfrak{a}_{\vec{\ell}, \lambda\,}]f \|_{L^p(\mathbb{R}^{d+1})} \leq C \lambda^{-\frac2p+ \varepsilon} 2^{c(\ell_2 +\cdots \ell_k)} \|f\|_{L^p(\mathbb{R}^{d})}.
	\]
\end{proposition}

To prove Proposition~\ref{prop_smoothing_reduc}, we consider two parts, $A_\Gamma[\mathfrak{a}_{\vec{\ell}, \lambda\,}\eta_\lambda ] $ and $A_\Gamma[\mathfrak{a}_{\vec{\ell}, \lambda\,}(1-\eta_\lambda) ]$, where
\begin{align}
\label{def:etalam}
	\eta_\lambda(y, \xi,\tau) 
	= \beta_0 \left(  \lambda^{-\varepsilon_1} ( \Gamma(y), 1)\cdot (\xi, \tau) \right) 
	\times \prod_{i=1}^{d-1} \beta_0 \left( \lambda^{-\frac12-\varepsilon_1} \partial_i \Gamma(y) \cdot  \xi\right),
\end{align}
for sufficiently small $\varepsilon_1$, to be chosen later. 
Note that $\mathfrak{D}^m\gamma$ is always larger than $\lambda^{-\frac12}$ and $\varepsilon_1$ can be chosen depending only on  $\varepsilon$ and $d$ (see \eqref{epsilon1} below).  
Then we have
\[
	A_\Gamma[\mathfrak{a}_{\vec{\ell}, \lambda\,}] = A_\Gamma[\mathfrak{a}_{\vec{\ell}, \lambda\,}\eta_\lambda ] + A_\Gamma[ \mathfrak{a}_{\vec{\ell}, \lambda\,} (1-\eta_\lambda )].
\]
We denote symbols by
\[
	\mathfrak{a}_{\vec{\ell}, \lambda}^{\text{in}} = \mathfrak{a}_{\vec{\ell}, \lambda\,}\eta_\lambda ,\quad 
	\mathfrak{a}_{\vec{\ell}, \lambda}^{\text{out}} = \mathfrak{a}_{\vec{\ell}, \lambda\,} (1-\eta_\lambda ),
\]
and the corresponding multipliers by
\begin{align*}
	\mathfrak{m}_{ \vec{\ell}, \lambda}^* (\xi, \tau) = &\iint \mathrm{e}^{-i s (\Gamma(y)\cdot\xi + \tau)} \mathfrak{a}_{\vec{\ell}, \lambda}^{*} ~\mathrm{d}y\mathrm{d}s,\quad * = \text{in, out}.
\end{align*}
Proposition \ref{prop_smoothing_reduc} directly follows from Lemma \ref{lem_out} and Lemma \ref{lem_decoup} below, which treat the \emph{out}-symbol and the \emph{in}-symbol, respectively.
For the \emph{out}-symbol, we use the lower bounds of $(\Gamma,1)\cdot(\xi, \tau)$ and $\partial_j \Gamma \cdot \xi$, which yield oscillatory integral estimates for the symbol.
For the \emph{in}-symbol, we will apply the smoothing estimates developed in \cite{Oh2025}, which exploits the geometry of localized frequency variables.

\begin{lemma}\label{lem_out}
    Let $p\geq2$ and $2^{\ell_m} \geq \lambda^{-1/2}$ for all $m=2, \dots, k$. Then, for some $c>0$, it holds that
    \[
        \left\| A_\Gamma[\mathfrak{a}_{\vec{\ell}, \lambda}^{\text{out}}] f \right\|_{L^p(\mathbb{R}^{d+1})}
        \lesssim \lambda^{-\frac{2}{p}} 2^{c(\ell_2+\cdots+\ell_k)}\|f\|_{L^p(\mathbb{R}^{d})}.
    \]
\end{lemma}

\begin{lemma}\label{lem_decoup}
    Let $p>4(k-1)$ and $2^{\ell_m} \geq \lambda^{-1/2}$ for all $m=2, \dots, k$. Then there are $C = C(d), c = c(k,p)>0$ such that
    \begin{equation}
        \nonumber
        \left\| A_\Gamma[\mathfrak{a}_{\vec{\ell}, \lambda}^{\text{in}}] f \right\|_{L^p(\mathbb{R}^{d+1})}
        \lesssim \lambda^{- \frac{2}{p}+C\varepsilon_1} 2^{c(\ell_2+\cdots+\ell_k)}\|f\|_{L^p(\mathbb{R}^{d})}.
    \end{equation}
\end{lemma}
Note that $\varepsilon_1$ is the parameter appearing in \eqref{def:etalam}. Thus Proposition \ref{prop_smoothing_reduc} follows by Lemmas \ref{lem_out} and \ref{lem_decoup} after choosing 
\begin{equation}
    \label{epsilon1}
 \varepsilon_1=\varepsilon/C(d).   
\end{equation}

\subsection{Proof of Lemma~\ref{lem_out}}

The proof of Lemma~\ref{lem_out} is completed by interpolation between $L^2$ and $L^\infty$ estimates, which follow from $L^\infty$-bound for $\mathfrak{m}_{ \vec{\ell}, \lambda}^{\text{out}} $ and $L^1$-estimates for $(\mathfrak{m}_{ \vec{\ell}, \lambda}^{\text{out}} )^{\vee}$, respectively.
That is, we show
\begin{align}
        	\| A_\Gamma [\mathfrak{a}_{\vec{\ell}, \lambda}^{\text{out}} ] f \|_{L^2(\mathbb{R}^{d+1})}\lesssim &\lambda^{-1} 2^{c(\ell_2 + \cdots + \ell_k)}\|f\|_{L^2(\mathbb{R}^{d})},\label{L2}\\
            \| A_\Gamma [\mathfrak{a}_{\vec{\ell}, \lambda}^{\text{out}} ] f \|_{L^\infty(\mathbb{R}^{d+1})}\lesssim  &2^{c(\ell_2 + \cdots + \ell_k)}\|f\|_{L^\infty(\mathbb{R}^{d})}.\label{Linfty}
\end{align}
We present proof of both estimates in the following subsections.

\subsubsection{Proof of \eqref{L2}: $L^\infty$-bounds for $\mathfrak{m}_{ \vec{\ell}, \lambda}^{\text{out}} $} 

We decompose $\mathfrak{a}_{\vec{\ell},\lambda}^{\text{out}}$ into two parts, distinguished by whether $| \Gamma(y)\cdot \xi +\tau | > \lambda^{\varepsilon_1}$ or $|\partial_i \Gamma(y) \cdot \xi | > \lambda^{\frac12+\varepsilon_1}$ for some $1\le i\le d-1$. For this purpose, define
\[
\begin{aligned}
    \mathfrak{a}_{\vec{\ell}, \lambda}^{\text{out}, 1}&=\mathfrak{a}_{\vec{\ell},\lambda}\cdot\big(1-\beta_0(\lambda^{-\epsilon_1}(\Gamma(y),1)\cdot(\xi,\tau))\big),\\
    	\mathfrak{a}_{ \vec{\ell}, \lambda}^{\text{out},2} &= \mathfrak{a}_{\vec\ell,\lambda}\cdot\big(\beta_0 \left( \lambda^{-\varepsilon_1} (\Gamma(y), 1)\cdot (\xi, \tau) \right)-\eta_\lambda(y,\xi,\tau)\big).
\end{aligned}
\]
Note that $\mathfrak{a}_{ \vec{\ell}, \lambda}^{\text{out}}=\mathfrak{a}_{ \vec{\ell}, \lambda}^{\text{out}, 1}+\mathfrak{a}_{ \vec{\ell}, \lambda}^{\text{out}, 2}$. 
Then \eqref{L2} is reduced to show
\begin{equation}
    \label{ineq:l2out1,2}
    	\| A_\Gamma [\mathfrak{a}_{\vec{\ell}, \lambda}^{\text{out},j} ] f \|_{L^2(\mathbb{R}^{d+1})}\lesssim \lambda^{-1} 2^{c(\ell_2 + \cdots + \ell_k)}\|f\|_{L^2(\mathbb{R}^{d})},\quad j=1,2.
\end{equation}

We first consider the symbol $\mathfrak{a}_{ \vec{\ell}, \lambda}^{\text{out}, 1}$ and denote by $\mathfrak{m}_{ \vec{\ell},\lambda}^{\text{out},1}$ the corresponding multiplier, defined similarly as before.
Since $| \Gamma(y)\cdot \xi +\tau | > \lambda^{\varepsilon_1}$ on $\supp (\mathfrak a_{\vl,\lambda}^{\text{out},1})$ and $s$-derivatives of $\mathfrak a_{\vl,\lambda}^{\text{out},1}$ are bounded by a constant, one could apply the integration by part in $s$-variable to show
\[
	| \mathfrak{m}_{ \vec{\ell}, \lambda}^{\text{out},1} | \lesssim \lambda^{-\varepsilon_1 N} \Big| \supp_{(y,s)}(\mathfrak{a}_{\vec{\ell}, \lambda}^{\text{out},1} ) \Big|,
\] 
where we denote $\supp_{(y,s)}\mathfrak a=\cup_{\xi,\tau}\supp(\mathfrak a(\cdot,\cdot,\xi,\tau))$. 
 By Plancherel's theorem, it follows that
\begin{align*}
	\| A_\Gamma [\mathfrak{a}_{\vec{\ell}, \lambda}^{\text{out},1} ] f \|_{L^2(\mathbb{R}^{d+1})}
	\lesssim &\Big(\int_{|\tau|\lesssim \lambda} \lambda^{-2 \varepsilon_1 N} \Big| \supp_{(y,s)}(\mathfrak{a}_{\vec{\ell}, \lambda}^{\text{out},1} ) \Big|^2~\mathrm{d}\tau\Big)^{1/2} ~ \| f\|_{L^2(\mathbb{R}^{d})}\\
	\lesssim & \lambda^{-\varepsilon_1 N +\frac12} \Big| \supp_{(y,s)}(\mathfrak{a}_{\vec{\ell}, \lambda}^{\text{out},1} ) \Big| \, \|f\|_{L^2(\mathbb{R}^{d})}.
\end{align*}
After choosing sufficiently large $N$ such that $\varepsilon_1 N > 3/2$,  the support size estimate \eqref{ineq:suppsize} gives that
\[
	 \lambda^{-\varepsilon_1 N +\frac12} \Big| \supp_{(y,s)}(\mathfrak{a}_{\vec{\ell}, \lambda}^{\text{out},1} ) \Big| \leq \lambda^{-1} \delta_{\vec{\ell}}.
\]
Thus, the desired bound \eqref{ineq:l2out1,2} with $j=1$ follows by \eqref{sscale_sumf}.

Now we focus on the symbol $\mathfrak{a}_{ \vec{\ell}, \lambda}^{\text{out},2}$. To verify \eqref{ineq:l2out1,2} with $j=2$,
 we observe that on $\supp 	(\mathfrak{a}_{\vec{\ell}, \lambda}^{\text{out},2}) $,
\begin{align}\label{ineq_out2}
	| \Gamma(y)\cdot \xi +\tau | \leq 2\lambda^{\varepsilon_1}\quad \text{and} \quad 
	|\partial_i \Gamma(y) \cdot \xi | > \lambda^{\frac12+\varepsilon_1}
\end{align}
for some $1\leq i\leq d-1$. 
We define a differential operator  $\mathcal{D}$ by
\[
	\mathcal{D} = \frac{\nabla_y(\Gamma(y)\cdot \xi) \cdot \nabla_y}{|\nabla_y (\Gamma(y)\cdot \xi)|^{2}}.
\]
After using integration by parts $N$ times, we obtain
\begin{equation}
\nonumber
\begin{aligned}
	\Big| \iint \mathrm{e}^{-is(\Gamma(y)\cdot\xi + \tau)} (\mathcal{D}^*)^N (\mathfrak{a}_{\vec{\ell}, \lambda}^{\text{out},2})~\mathrm{d}y\mathrm{d}s \Big| 
	\leq \sup_{y,s} |(\mathcal{D}^*)^N (\mathfrak{a}_{\vec{\ell}, \lambda}^{\text{out},2})| \times |\supp_{(y,s)}(\mathfrak{a}_{\vec{\ell}, \lambda}^{\text{out},2})|.
\end{aligned}
\end{equation}
Using the observation $|\supp_{(y,s)}(\mathfrak{a}_{\vec{\ell}, \lambda}^{\text{out},2})| \lesssim \delta_{\vec{\ell}\,}$ together with \eqref{sscale_sumf}, to prove \eqref{ineq:l2out1,2} in the case $j=2$, it suffices to show that 
\begin{align}
    \label{ineq:allout2}
    |(\mathcal{D}^*)^N (\mathfrak{a}_{\vec{\ell}, \lambda}^{\text{out},2})|\lesssim \lambda^{-1}.
\end{align}

To estimate $(\mathcal{D}^*)^N (\mathfrak{a}_{\vec{\ell}, \lambda}^{\text{out},2})$, 
recall that \eqref{def:allambda} and \eqref{def:etalam}, which gives
\[
\begin{aligned}
    \mathfrak{a}_{\vec{\ell},\lambda}^{\text{out},2}=\mathfrak b(s,\xi,\tau)\psi_{\vec{\ell}\,}(y)\beta_0 \left( \lambda^{-\varepsilon_1} (\Gamma(y), 1)\cdot (\xi, \tau) \right)\big(1-\prod_{i=1}^{d-1} \beta_0( \lambda^{-\frac12-\varepsilon_1} \partial_i \Gamma(y) \cdot  \xi)\big),
\end{aligned}
\]
where $\mathfrak b(s,\xi,\tau)=\zeta(s) \beta_1(\lambda^{-1}|\xi|) \beta_0((C_\circ\lambda)^{-1}\tau)$.
We need to check the derivatives of each component of $\mathfrak{a}_{\vec{\ell}, \lambda}^{\text{out},2} $ and the term $\nabla_y (\Gamma(y)\cdot \xi)\times |\nabla_y (\Gamma(y)\cdot \xi)|^{-2}$ in the definition of $\mathcal D$, with respect to $y$.
First we have
\begin{align}
\nonumber
	| \partial_y^\alpha \psi_{\vec{\ell}\,} | \lesssim \delta_{\vec{\ell}\,}^{-|\alpha|}.
\end{align}
Secondly, 
we observe that on $\supp (\mathfrak{a}_{\vec{\ell}, \lambda}^{\text{out},2})$,
\begin{align}
	|\partial_y^{\alpha} (\lambda^{-\varepsilon_1} (\Gamma, 1)\cdot (\xi, \tau))| 
	\lesssim & (\lambda^{-\varepsilon_1}| \nabla_y(\Gamma\cdot\xi)| )^{|\alpha|},\quad |\alpha|\ge0,\label{250321_1559}\\
	|\partial_y^{\alpha} (\lambda^{-\frac12 -\varepsilon_1} \partial_i\Gamma\cdot\xi)|
	\lesssim &  (\lambda^{-\varepsilon_1} |\nabla_y (\Gamma\cdot\xi )|)^{|\alpha|},\quad |\alpha|\ge1.\label{250321_1603}
\end{align}
Indeed, \eqref{250321_1559} is trivial for $|\alpha|=0,1$. For $|\alpha|\geq2$, we get
\begin{equation}
\label{ineq:partialalphagamma}
    	|\partial_y^{\alpha} ( (\Gamma, 1)\cdot (\xi, \tau))| 
	\lesssim  \lambda
	\le (\lambda^{-\varepsilon_1}| \nabla_y(\Gamma\cdot\xi)| )^{|\alpha|},
\end{equation}
due to $|(\xi, \tau)|\lesssim  \lambda$, \eqref{upperbounds:gammma}, and $|\nabla_y (\Gamma \cdot \xi) | > \lambda^{\frac12+\varepsilon_1}$.
Analogously to \eqref{250321_1559}, \eqref{250321_1603} follows directly from \eqref{ineq:partialalphagamma}.
Thus we have 
\begin{align}\label{int by part 1}
\begin{split}
	|\partial_y^\alpha \mathfrak{a}_{\vec{\ell}, \lambda}^{\text{out},2} | 
	\lesssim  &\sum_{\alpha = \alpha_1 + \alpha_2} \delta_{\vec{\ell}}^{-|\alpha_1|} \times (\lambda^{-\varepsilon_1}| \nabla_y(\Gamma\cdot\xi)| )^{|\alpha_2|}.
\end{split}
\end{align}
On the other hand, by \eqref{ineq:partialalphagamma}, we also obtain for $|\alpha|\ge0$ and $i=1,\cdots d-1$,
\begin{align}\label{int by part 2}
	\Big| \partial_y^\alpha \left( \frac{\partial_i \Gamma \cdot \xi}{|\nabla_y \Gamma\cdot\xi|^{2}}   \right) \Big|
	\lesssim \frac{( \lambda^{-\varepsilon_1} |\nabla_y (\Gamma\cdot\xi)|)^{|\alpha|}}{|\nabla_y (\Gamma\cdot\xi)|}.
\end{align}

Now, one could compute derivatives of $\mathfrak{a}_{\vec{\ell}, \lambda}^\text{out,2}$,
\begin{align*}
	|(\mathcal{D}^*)^N (\mathfrak{a}_{\vec{\ell}, \lambda}^{\text{out},2})|
	\lesssim \sum_{|\alpha_1| + \cdots |\alpha_{N+1}| = N}
	|\partial_y^{\alpha_{N+1}} \mathfrak{a}_{\vec{\ell}, \lambda}^{\text{out},2} | \prod_{i=1}^N\Big| \partial_y^{\alpha_i} \left( \frac{\nabla_y \Gamma \cdot \xi}{|\nabla_y \Gamma\cdot\xi|^{2}}   \right) \Big|.
\end{align*}
From \eqref{int by part 1} and \eqref{int by part 2}, the summand in the right-hand side is less than a constant multiple of 
\begin{align*}
	&\sum_{\beta_1 + \beta_2 = \alpha_{N+1}} \delta_{\vec{\ell}}^{-|\beta_1|} \times (\lambda^{-\varepsilon_1}| \nabla_y(\Gamma\cdot\xi)| )^{|\beta_2|}\prod_{i=1}^N\Big| \partial_y^{\alpha_i} \left( \frac{\nabla_y \Gamma \cdot \xi}{|\nabla_y \Gamma\cdot\xi|^{2}}  \right) \Big|\\
	\lesssim &\sum_{\beta_1 + \beta_2  = \alpha_{N+1}} \delta_{\vec{\ell}}^{-|\beta_1|} \times (\lambda^{-\varepsilon_1}| \nabla_y(\Gamma\cdot\xi)| )^{|\beta_2|}  \frac{( \lambda^{-\varepsilon_1} |\nabla_y (\Gamma\cdot\xi)|)^{|\alpha_1|+\cdots+|\alpha_N|} } { |\nabla_y (\Gamma\cdot\xi)|^{N} }\\
	=&\sum_{\beta_1 + \beta_2  = \alpha_{N+1}} \delta_{\vec{\ell}}^{-|\beta_1|} \lambda^{-\varepsilon_1(N-|\beta_1|)}  |\nabla_y (\Gamma\cdot\xi)|^{-|\beta_1|}.
\end{align*}
Making use of \eqref{ineq_out2}, it follows that
\begin{align}
\nonumber
	 \delta_{\vec{\ell}}^{-|\beta_1|} \lambda^{-\varepsilon_1(N-|\beta_1|)}  |\nabla_y (\Gamma\cdot\xi)|^{-|\beta_1|}
	 \lesssim (\lambda^{\frac12} \delta_{\vec{\ell}} )^{-|\beta_1|} \lambda^{-\varepsilon_1 N} ,
\end{align}
which implies
\[
|(\mathcal{D}^*)^N (\mathfrak{a}_{\vec{\ell}, \lambda}^{\text{out},2})|\lesssim \sum_{|\beta_1|\leq N}( \lambda^{\frac12}\delta_{\vec{\ell}} )^{-|\beta_1|} \lambda^{-\varepsilon_1 N}.
\]
From $\lambda^{-1/2}\leq 2^{\ell_m}$ for all $m=2,\dots, k$, we have $\delta_{\vec \ell \,}\gtrsim  \lambda^{-1/2}$. Thus \eqref{ineq:allout2} follows after choosing $N$ sufficiently large.

\subsubsection{Proof of \eqref{Linfty}: Kernel estimates for $\mathfrak{a}_{\vec{\ell}, \lambda}^{\text{*}}$}
So far, the bounds for the symbols only produce $L^2$ estimates. 
The desired $L^p$ bounds follow by $L^\infty$ bounds based on the kernel estimates.
Since we also need the kernel estimate for the \emph{in}-symbol later, we give $L^1$ estimates of kernels for both the \emph{in}- and \emph{out}-symbols.  Define $K_{ \vec{\ell}, \lambda}^{\text{*}}$ for $*=$ in, out, as 
\begin{align*}
	K_{ \vec{\ell}, \lambda}^{\text{*}}(x,t) 
	= &\int \mathcal{K}[\mathfrak a_{\vec{\ell}, \lambda}^{\text{*}}] (x, t; y)~\mathrm{d}y,\\
    \mathcal{K}[\mathfrak a_{\vec{\ell}, \lambda}^{\text{*}}] (x, t; y)
    \coloneq &\iiint \mathrm{e}^{2\pi i \left( (x- s\Gamma(y))\cdot\xi + (t-s)\tau \right)} \mathfrak{a}_{\vec{\ell}, \lambda}^{\text{*}}(y,s,\xi,\tau)~\mathrm{d}s\mathrm{d}\xi\mathrm{d}\tau.
\end{align*}
Then $K_{ \vec{\ell}, \lambda}^{\text{*}}$ is the inverse Fourier transform of the multiplier $\mathfrak{m}_{ \vec{\ell}, \lambda}^{\text{*}}$ and satisfies the following lemma.

\begin{lemma}\label{lem_kernel}
	For $*=$ in,out, we have
	\[
		\int_{\mathbb{R}^d} |K_{ \vec{\ell}, \lambda}^{\text{*}}|~\mathrm{d}x \lesssim |\supp (\psi_{\vec{\ell}}\,)|.
	\]
\end{lemma}

By Lemma~\ref{lem_kernel} and Young's inequality, one easily obtains $L^\infty$ bounds for $*=$in, out, 
\begin{equation}
    \label{est:linfty}
    		\| A_\Gamma[\mathfrak{a}^{\text{*}}_{\vec{\ell}, \lambda\,}]f \|_{L^\infty(\mathbb{R}^{d+1})}\lesssim |\supp(\psi_\vl)|\|f\|_{L^\infty(\mathbb{R}^{d})} \lesssim   \delta_\vl \,\|f\|_{L^\infty(\mathbb{R}^{d})}.
\end{equation}
Then Lemma \ref{lem_out} follows directly from
interpolation between $L^2$ bounds \eqref{ineq:l2out1,2} and $L^\infty$ bounds \eqref{est:linfty} combined with \eqref{sscale_sumf}.
\begin{proof}[Proof of Lemma \ref{lem_kernel}]
Note that $\mathfrak{a}_{\vec{\ell}, \lambda}^{\text{in}}=\mathfrak{a}_{\vec{\ell}, \lambda} \eta_\lambda$, $\mathfrak{a}_{\vec{\ell}, \lambda}^{\text{out}} =  \mathfrak{a}_{\vec{\ell}, \lambda} -\mathfrak{a}_{\vec{\ell}, \lambda} \eta_\lambda$, respectively. To prove Lemma~\ref{lem_kernel}, it suffices to show
\[
    \iint \Big|\mathcal{K}[\mathfrak a_{\vec{\ell}, \lambda}] (x, t; y)\Big|~\mathrm{d}y\mathrm dx,\, \iint \Big|\mathcal{K}[\mathfrak a_{\vec{\ell}, \lambda}\eta_\lambda] (x, t; y)\Big|~\mathrm{d}y\mathrm dx \lesssim |\supp (\psi_{\vec{\ell}}\,)|.
\]
For the integral of the kernel $\mathcal{K}[\mathfrak{a}_{\vec{\ell}, \lambda}](x,t;y)$, one has
\begin{align*}
	&\int_{\bR^d}  \int_{\bR^{d-1}} \Big| \mathcal{K}[\mathfrak{a}_{\vec{\ell}, \lambda}](x,t;y)\Big|~\mathrm{d}y\,\mathrm{d}x\\
	\lesssim &\int \iint  \Big| \psi_{\vec{\ell}\,}(y) \zeta(s) \lambda^d \widehat{\beta_1(|\cdot|)}\left(\lambda(x - s\Gamma(y))\right) \lambda \widehat{\beta}_0(C_\circ\lambda(t-s)) \Big|~\mathrm{d}y\mathrm{d}s\,\mathrm{d}x\\
	\lesssim &\int_{\bR^{d-1}} | \psi_{\vec{\ell}\,}(y) |~\mathrm{d}y.
\end{align*}
Note that integrals with respect to $x$ and $s$ are harmless, since $\beta_0, \beta_1$ are $C_c^\infty$ functions.

For $\iint  | \mathcal{K}[\mathfrak a_{\vec{\ell}, \lambda} \eta_\lambda ]|~\mathrm{d}y\mathrm{d}x$, we recall the following lemma of the last author \cite{Oh2025}.
\begin{lemma}[\cite{Oh2025}, Lemma~2.4.]\label{lem_kbound_Oh}
	Let $g : \mathbb{R}^{d-1} \times \mathbb{R} \to \mathbb{R}$ be a function such that $|g|\leq 1$, and set $\mathfrak{a}(y, s, \xi, \tau) = g(y, s) \beta_1(\lambda^{-1}|\xi|)\beta_0((C_\circ\lambda)^{-1}\tau) \eta_\lambda(y, \xi, \tau)$.
	Then we have
	\[
		\int \left| \mathcal{K}[\mathfrak{a}](x, t; y) \right|~\mathrm{d}x \lesssim 1.
	\]
\end{lemma}
\noindent By Lemma~\ref{lem_kbound_Oh} with $g(y,s)=\zeta(s)\psi_\vl(y)$, we have
\[
    \iint  | \mathcal{K}[\mathfrak a_{\vec{\ell}, \lambda} \eta_\lambda ]|~\mathrm{d}y\mathrm{d}x \lesssim \int_{\bR^{d-1}} | \psi_{\vec{\ell}\,}(y) |~\mathrm{d}y
\]
as desired.    
\end{proof}

\subsection{Proof of Lemma~\ref{lem_decoup}}

For the \emph{in}-symbol $\mathfrak{a}_{\vec{\ell}, \lambda}^{\text{in}}$, we introduce smoothing estimates of $A_\Gamma[\mathfrak{a}_{\vec{\ell}, \lambda}^{\text{in}}]$ for some $p$ in terms of $\lambda$ with decay in $\ell_2, \dots, \ell_k$.
We note that the decomposition of \emph{in}-symbol allows us to consider frequency variables $(\xi, \tau)$ are restricted in a conic surface in the sense of local geometry.
Thus, we can directly apply certain smoothing estimates arising from decoupling inequality due to such local geometry.

    Recall that we need to show  
    \begin{equation}
        \label{ineq:lem;decoup}
        \left\| A_\Gamma[\mathfrak{a}_{\vec{\ell}, \lambda}^{\text{in}}] f \right\|_{L^p(\mathbb{R}^{d+1})}
        \lesssim \lambda^{- \frac{2}{p}+C\varepsilon_1} 2^{c(\ell_2+\cdots+\ell_k)}\|f\|_{L^p(\mathbb{R}^{d})}
    \end{equation}
    for $p>4(k-1)$ and $2^{\ell_m} \geq \lambda^{-1/2}$ for all $m=2, \dots, k$ with $C = C(d), c = c(k,p)>0$.
    In \cite[Proposition~3.2]{Oh2025}, it was proved that there exists a constant $C$ such that \eqref{ineq:lem;decoup} holds with $c=0$ and $p>\max\{4(k-1),6\}$. Since we assume that $k\ge3$, we have $\max\{4(k-1),6\}=4(k-1)$. On the other hand, by  \eqref{sscale_sumf} and \eqref{est:linfty}, we have 
        \[
        \left\| A_\Gamma[\mathfrak{a}_{\vec{\ell}, \lambda}^{\text{in}}] f \right\|_{L^\infty(\mathbb{R}^{d+1})}
        \lesssim 2^{c(\ell_2+\cdots+\ell_k)}\|f\|_{L^\infty(\mathbb{R}^{d})},
    \]
    for some $c>0$.
    Thus, interpolation between this and \eqref{ineq:lem;decoup} with $c=0$ gives the desired estimate.

\begin{remark}
A symbol appeared in \cite[Proposition~3.2]{Oh2025} is slightly different from $\mathfrak{a}_{\vec{\ell}, \lambda}^{\text{in}}$.
Precisely, we decompose the support of $\psi$ so that one has $\mathfrak D^m \gamma (y) \sim 2^{\ell_m}$, while the decomposition of \cite{Oh2025} is associated with the scale $\mathfrak{D}^m\gamma(y) \sim \lambda^{-1/2} 2^{\ell_m}$. The difference would be easily overcome by taking $\lambda^{1/2} 2^{\ell_m}$ instead of $2^{\ell_m}$ in the proof of \cite[Proposition~3.2]{Oh2025}, so we obtain \eqref{ineq:lem;decoup} with $c=0$.
\end{remark}

\section{Maximal estimates, sublevel sets, and Fourier decay}\label{sec_fdecay}

In this section, we prove Theorem~\ref{thm_Fdecay}. To this end, we first study the interplay between $L^p$ bounds of maximal functions and sublevel set estimates.
Moreover, we extend the classical sublevel set estimates derived from a result of Iosevich--Sawyer \cite{IS1}  so that the estimates remain valid uniformly under translation by arbitrary affine functions. Secondly, we adapt the stationary set method of Basu--Guo--Zhang--Zorin-Kranich \cite{BGZZ} to estimate the Fourier transform of measures on analytic hypersurfaces. 
Combining the stationary set method and the extended sublevel set estimates, we obtain the desired estimates for oscillatory integrals, which are the Fourier decay of analytic surfaces.

To proceed as mentioned above, we first give a simple reduction of Theorem \ref{thm_Fdecay}. By compactness of the hypersurface $\Gamma$, the problem reduces to the case when $\Gamma$ is given as the graph of an analytic function $\gamma:U\rightarrow \bR$ over an open set $U\subset \bR^{d-1}$. 
In this setting, the associated maximal operator,
\[
    M_\Gamma[\mathds1_{U}]f(x) := \sup_{t>0} \Big| \int_{U} f\left(x' - ty, x_d - t\gamma(y) \right) \,\mathrm dy \Big|,\quad x= (x', x_d)\in\mathbb R^d,
\]
is bounded on $L^p(\bR^d)$ for $p>p_{cr}$, by the hypothesis of Theorem \ref{thm_Fdecay}. Under this assumption, to prove Theorem \ref{thm_Fdecay}, it suffices to show that 
\[
|\mathcal F[\mu[\psi]](\xi)|\lesssim |\xi|^{-1/p}, 
\]
whenever $\psi\in C_c^\infty(U)$ vanishes on an open set containing all points $y$ for which $(y,\gamma(y))$ is non-transversal (see \eqref{def:mupsi} for the definition of $\mu[\psi]$). This condition on $\psi$ can be equivalently expressed in terms of the quantity,
    $$c_\psi:=\inf_{v\in \supp (\psi)}|\nabla\gamma(v)\cdot v-\gamma(v)|>0.$$ 
Indeed, for each $v$, the quantity
\begin{align}
\nonumber
    c_v := &-\gamma(v) +\nabla\gamma(v)\cdot v,
\end{align}
 vanishes precisely when $(v, \gamma(v))$ is a non-transversal point, and is nonzero otherwise. Thus, the requirement that $\psi$ vanishes near non-transversal part is equivalent to the assumption $c_\psi>0$.
Consequently, Theorem \ref{thm_Fdecay} is deduced from the following theorem.

\begin{theorem}\label{thm_Fdecay2}
    Let $d\ge2$, $p_{cr}\ge 1$, $\gamma$ be an analytic function on an open set $U\subset \mathbb B_1^{d-1}(0)$, and $\psi\in C_c^\infty(U)$ satisfy $c_\psi>0$.
    Suppose that $M_\Gamma[\mathds1_{U}]$ is bounded on $L^p(\bR^{d})$ for $p>p_{cr}$. Then, for any $p>p_{cr}$, we have
    \[
    |\mathcal F[\mu[\psi]](\xi)|\lesssim (c_\psi|\xi|)^{-\frac 1p}.
    \]
\end{theorem}
In the remainder of this section, we devote ourselves to proving Theorem \ref{thm_Fdecay2}.
To this end, we first establish certain sublevel set estimates from maximal estimates, then obtain the desired Fourier decay by making use of the sublevel set estimates.

\subsection{From maximal estimates to sublevel set estimates}

Iosevich--Sawyer \cite[Theorem 2]{IS1} showed that  $M_{\Gamma}[\mathds1_{U}]$ is bounded on $L^p $ only if $\dist(\cdot, H)^{-1} \in L_{\text{loc}}^{1/p}(\Gamma, \mathrm d\mu) $ for all  tangent hyperplane $H$ of $\Gamma$ at $(y_H, \gamma(y_H))$ not passing through the origin. 
Suppose $\chi$ is a  cut-off function around $y_H$.
Then, the local integrability of $\dist(\,\cdot, H)^{-1/p}$ is equivalent to the following inequality: 
\begin{align}\label{ineq_crtcl}
    \int |\gamma(y) - \gamma(y_H) - \nabla\gamma(y_H)\cdot(y-y_H)|^{-1/p} \chi(y) \,\mathrm d y \le C.
\end{align}
Indeed, consider a point $y_H \in \mathbb B_1^{d-1}(0)$ and a hyperplane $H$ tangent to $\Gamma$ at $(y_H,\gamma(y_H))$ so that $(y, y_d) \in H$ satisfies
    \[
        y_{d} = \gamma(y_H) + \nabla\gamma(y_H)\cdot(y-y_H).
    \]
    For a point $(y, \gamma(y))$ on $\Gamma$, the function $\dist((y, \gamma(y)), H)$, is given by
    \[
        | (y - y_H, \gamma(y) - \gamma(y_H)) \cdot \boldsymbol{n}_H |,
    \]
    where $\boldsymbol{n}_H$ denotes a unit normal vector of $H$.
    Observe that $\boldsymbol{n}_H$ is obtained by normalizing $(\nabla\gamma(y_H), -1)$, and one has
    \[
        |(y - y_H, \gamma(y) - \gamma(y_H))\cdot (\nabla\gamma(y_H), -1) | = |\gamma(y) - \gamma(y_H) - \nabla\gamma(y_H)\cdot(y-y_H)|.
    \]
    Thus, $\dist( (y, \gamma(y)), H)$ and $|\gamma(y) - \gamma(y_H) - \nabla\gamma(y_H)\cdot(y-y_H)|$ are equivalent up to a constant depending on $|(\nabla\gamma(y_H), -1)|$.
Moreover, \eqref{ineq_crtcl} implies the following sublevel set estimate by Chebyshev's inequality: 
\begin{align}\label{ineq_sls}
    |\{y \in \supp(\chi) : |\gamma(y)-\gamma(y_H)-\nabla\gamma(y_H)\cdot(y-y_H)|\le \lambda^{-1} \}| \le C\lambda^{-1/p}.
\end{align}
The converse is generally not true. (See \cite[Section~1.1]{CGM2013}.)

To obtain the Fourier decay of ${\mu[\psi]}$, we require decay estimates for oscillatory integrals with phases of the form $\lambda(\gamma(y)+u\cdot y)$ uniformly with respect to $u$.
The main goal of this subsection is to extend \eqref{ineq_crtcl} and \eqref{ineq_sls} to \eqref{ineq_crtcl_pert} and \eqref{ineq:sublevelset} below, which are closely related to the uniform estimates.
To this end, we use the observation that the maximal average has an
$L^p$ operator norm invariant under the linear transform $T:(y, y_{d})\mapsto (y, y_{d}+u\cdot y)$.
 More precisely, one has $\|f\circ T\|_p = \|f\|_p$ and
\[
    \| M_\Gamma[\mathds 1_U] (f\circ T)\|_p = \| M_{\Gamma_u}[\mathds 1_U]f \circ T \|_p = \|M_{\Gamma_u}[\mathds 1_U]f\|_p,
\]
where $\Gamma_u$ denotes the hypersurface $\{(y, \gamma(y) + u\cdot y):y\in U\}$.
Consequently, $M_{\Gamma_u}[\mathds 1_U]$ is bounded on $L^p(\mathbb R^d)$ with the same operator norm as $M_\Gamma[\mathds 1_U]$. Heuristically, this invariance reflects the fact that the map $T$ preserves the underlying transversality structure, which affects the $L^p$ boundedness of $M_\Gamma[\mathds 1_U]$.

Utilizing the above observation, we prove the following theorem, which provides uniform estimates of \eqref{ineq_crtcl} over affine shifts.
\begin{theorem}\label{thm_lin_pert}
    Let $p_0>1$ and  $U\subset \mathbb B_1^{d-1}(0)$ be an open set. 
    Suppose that $M_\Gamma[\mathds1_U]$ is bounded on $L^{p_0}(\mathbb R^{d})$.
    Then, for any $p>p_0$ and $v\in U$, there is a neighborhood $U_v$ of $v$ such that for all $u\in \bR^{d-1}$,
    \begin{align}\label{ineq_crtcl_pert}
        \int_{U_v} |\gamma(y)-\gamma(v) - \nabla\gamma(v)\cdot(y-v) + u\cdot (y-v)|^{-1/p} \,\mathrm dy \le C |c_v|^{-1/p},
    \end{align} 
    where $C$ is an implicit constant independent of $u, v$.
\end{theorem}

We may apply Theorem \ref{thm_lin_pert} to obtain \eqref{ineq_crtcl_pert} uniformly, since $|c_v|\ge c_\psi$ for all $v\in \supp(\psi)$.
We note that Iosevich--Sawyer also showed \eqref{ineq_crtcl_pert}, but the dependence on affine terms has not been explicit. 
That is, Theorem~\ref{thm_lin_pert} gives a uniform upper bound with respect to the affine terms  $u\cdot(y-v)$, which is crucial in our proof of Theorem~\ref{thm_Fdecay2}.

\begin{proof}
Since \eqref{ineq_crtcl_pert} is trivial when $c_v=0$, we assume that $c_v\neq 0$.
The proof splits into two cases: $|u|\ge |c_v|/10$ or $|u|<|c_v|/10$. For the first case, let $U_v$ be a small ball centered at $v$ satisfying
\begin{equation}
      \label{ineq:nablagamma}
|\nabla\gamma(y)-\nabla\gamma(v)|\le|c_v|/100,      
\end{equation}
for all $y\in U_v$. The existence of such a ball  is guaranteed by Taylor's theorem. Then the gradient of $\gamma(y)-\gamma(v)-(\nabla\gamma(v)-u)\cdot(y-v)$  with respect to $y$ satisfies the lower bound
\begin{equation*}  
|\nabla\gamma(y)-\nabla\gamma(v)+u|\gtrsim |c_v|,
\end{equation*}
for $y\in U_v$. By Fubini's theorem and the lower bound on the gradient, we can show that 
\[
\int_{U_v} |\gamma(y)-\gamma(v) - \nabla\gamma(v)\cdot(y-v) + u\cdot (y-v)|^{-1/p} \,\mathrm dy\lesssim \int_{-1}^1 |c_vt|^{-1/p}\,\mathrm dt,
\]
which gives \eqref{ineq_crtcl_pert}. 
Indeed, by \eqref{ineq:nablagamma}, there exists a unit vector $\omega $ satisfying $|(\nabla\gamma(y)-\nabla\gamma(v)+u)\cdot \omega|\gtrsim |c_v|$ for all $y\in U_v$. After a suitable rotation, one may assume that 
\begin{equation}
    \label{ineq:e1directionalbound}
    |(\nabla\gamma(y)-\nabla\gamma(v)+u)\cdot e_1|\gtrsim |c_v|.
\end{equation}
 For fixed $y_2,\cdots, y_{d-1}$, since a set $\{y_1:y\in U_v\}$ is an interval, there exists $y_1^*$ satisfying
\[
\inf_{y_1:y\in U_v}|\gamma_v(y-v)+ u\cdot (y-v)|=|\gamma_v(y^*-v)+ u\cdot (y^*-v)|,
\]
    where $y^*=(y_1^*,y_2,\cdots,y_{d-1})$ and  $\gamma_v$ is given by
    \begin{align*}
        \gamma_v (y) \coloneq &\gamma(y+ v)-\gamma(v) - \nabla \gamma(v) \cdot y.
    \end{align*}
    Then, by \eqref{ineq:e1directionalbound}, we have
\[
\int_{y_1:y\in U_v} |\gamma_v(y-v) + u\cdot (y-v)|^{-1/p} \,\mathrm dy_1\lesssim \int_{-2}^2|c_v(y_1-y_1^*)|^{-1/p} \,\mathrm dy_1,
\]
which gives the desired inequality as a consequence of Fubini's theorem.

Now, we focus on the second case, $|u|<c_v/10$.
 For this case, we additionally assume that $U_v$ satisfies   $U_v\subset U$,
 \begin{align}\label{ineq_covering}
    \begin{split}
        |\gamma_v(y-v)| &\le |c_v|/100,\\
        |y-v| &\le 1/100,
    \end{split}
    \end{align}
    for all $y \in U_v$. Such a ball $U_v$ can be chosen by Taylor's theorem.
    After translation to the origin, one can consider $y$ being contained in $U_v':=-v+U_v$, which is a ball centered at the origin.
    That is to say, we consider the following maximal operator:
    \[
        \widetilde M_{u,v} f (x) = \sup_{t>0} \Big|\int_{\mathbb R^{d-1}} f\left(x' - t(y+v), x_d -t(\gamma_v(y) + u\cdot (y+ v) -c_{v})\right) \mathds1_{ U_v'}(y)\, \mathrm dy \Big|.
    \]
    Note that $\|\widetilde M_{u,v} f\|_{p_0}\le \|\widetilde M_{u,v} |f|\|_{p_0}\le \|M_\Gamma[\mathds1_U] |f|\|_{p_0}$. 
    Hence, the norm $\|\widetilde M_{u,v} f\|_{p_0}$  is bounded by a constant multiple of $\|f\|_{p_0}$ uniformly with respect to $u,v$. Then, we show that the LHS of \eqref{ineq_crtcl_pert} is bounded by $\|\widetilde M_{u,v} f\|_{p_0}$ for some specific $f$. Furthermore, with a suitable choice of $f$ and a careful analysis, we conclude that \eqref{ineq_crtcl_pert} indeed holds.

    To be precise, we take 
    \begin{align}
    \label{def:ff}
        f(x', x_d) = |x_d|^{-1/p} \mathds1_{\mathbb B_{1}^{d-1}(0)}(x') \mathds1_{[-1, 1]}(x_d/c_v).
    \end{align}
    Note that $\|f\|_{p_0}$ is controlled by a constant multiple of $|c_v|^{\frac1{p_0} - \frac1{p}}(1- \frac{p_0}{p})^{-1/{p_0}}$.
    Now, when $x_d/c_v$ is negative, choose 
    $$t = -x_d/(c_v-u\cdot v).$$ 
    Then it is positive since both $c_v-u\cdot v$ and $c_v$ have the same sign from $|u|<c_v/10$ and $|v|<1$.  By this observation, for such $x$, and with the choice of $t$, the maximal function $\widetilde M_{u,v} f(x)$ is bounded below by 
    \[
    \int_{U_v'}  (t|\gamma_v(y) + u\cdot y|)^{-1/p} \mathds1_{\mathbb B_1^{d-1}(0)}(x'-t(y+v)) \mathds1_{[-1, 1]}(-t (\gamma_v(y)+ u \cdot y)/c_v)\mathrm{d}y.
    \]
    Moreover,
    we can show that 
    \begin{align}
    \label{ineq:lbtildemgamma}
            \widetilde M_{u,v} f(x)\ge\mathds1_{\mathbb B_{1/2}^{d-1}(0)}(x') \mathds1_{[-10^{-1}, 0)}(x_d/c_v) \int_{U_v'} (t|\gamma_v(y) + u\cdot y|)^{-1/p} \,\mathrm dy.
        \end{align}
    Indeed, for  $x_d/c_v\in [-10^{-1},0)$, $y\in U_v'$, $v\in V\subset \mathbb B^{d-1}_1(0)$, and $u\in \mathbb B_{|c_v|/10}^{d-1}(0)$, 
    one has
    \begin{align}
    \label{ineq:x'}
        \mathds1_{\mathbb B_1^{d-1}(0)}(x'-t(y+v)) \ge &\mathds1_{\mathbb B_{1/2}^{d-1}(0)}(x'),\\
        \nonumber
        \mathds1_{[-1, 1]}(-t (\gamma_v(y)+ u \cdot y)/c_v) \ge & 1,
    \end{align}
 since $|x_d|/2\le |c_v|t\le 2|x_d|$ and \eqref{ineq_covering} hold.
   
    After taking $L^{p_0}$-norm on both sides of \eqref{ineq:lbtildemgamma},  it follows that
    \[
        \| \widetilde M_{u,v} f\|_{p_0} \gtrsim |c_v|^{1/{p_0}} (1- \frac{p_0}{p})^{-1/{p_0}} \int_{U_v'} |\gamma_v(y)+ u\cdot y|^{-1/p} \,\mathrm dy.
    \]
    By the $L^{p_0}$ boundedness of $M_\Gamma[\mathds1_U]$, we have
        \[
            \| \widetilde M_{u,v} f\|_{p_0} \le K \| f\|_{p_0} \lesssim  K |c_v|^{\frac1{p_0} - \frac1{p}}(1- \frac{p_0}{p})^{-1/{p_0}},
        \]
    where $K$ denotes the constant such that $\|M_\Gamma[\mathds1_U]\|_{p_0\to p_0} \le K$.
    Thus, we conclude that
    \begin{align}
    \nonumber
        \int_{U_v'} | \gamma_v(y)+u\cdot y |^{-1/p}  \,\mathrm dy \le C K |c_v|^{-\frac1{p}},
    \end{align}
    whenever $u \in \mathbb B_{|c_v|/10}^{d-1}(0)$, which yields \eqref{ineq_crtcl_pert}.
\end{proof}

Once we have \eqref{ineq_crtcl_pert}, one can directly obtain associated sublevel set estimates by applying Chebyshev's inequality. 
In what follows, we obtain a stronger version of these estimates that holds uniformly under arbitrary constant translations of the original function.
This translation-uniformity is essential for the results in the following section.

\begin{theorem}\label{thm_lin_pert2}
    Under the hypothesis of Theorem \ref{thm_lin_pert},
     for any $p>p_0$ and $v\in U$,  there is a neighborhood $U_v$ of $v$  such that for all $\lambda\ge1$, $u\in \bR^{d-1}$, and $\beta\in\bR$,
    \begin{equation}
        \label{ineq:sublevelset}
        |\{y\in U_v : |\gamma(y) + u\cdot y - \beta| \le \lambda^{-1}\}| \lesssim  |c_v|^{-1/p}\lambda^{-1/p}.
    \end{equation}
    Here, the implicit constant is independent of $\lambda, u, v, \beta$.
\end{theorem}
\begin{proof}
 As in the proof of Theorem \ref{thm_lin_pert}, we choose $U_v$ to be a ball centered at $v$ satisfying \eqref{ineq:nablagamma} and \eqref{ineq_covering}.
When $|u+\nabla\gamma(v)|\ge|c_v|/10$, the same argument as in the first part of the proof of Theorem \ref{thm_lin_pert} implies
\[
\int_{U_v} |\gamma(y) + u\cdot y-\beta|^{-1/p} \,\mathrm dy\lesssim |c_v|^{-1/p}.
\]
Thus, the desired follows by Chebyshev's inequality.

For the case $|u+\nabla\gamma(v)|<|c_v|/10$, we proceed as in the case $|u|<|c_v|/10$ of the proof of Theorem \ref{thm_lin_pert}, but with a different choice of $t$. To be precise, set $\widetilde u=u+\nabla\gamma(v)$ and $\widetilde \beta=\beta+c_v-\widetilde u\cdot v$, so that 
\begin{equation}
\label{identity:gammaub}
\gamma(y)+u\cdot y-\beta=\gamma_v(y-v)+\widetilde u\cdot (y-v)-\widetilde \beta.    
\end{equation}
On $U_v$, we observe that $|\gamma_v(y-v)+\widetilde u\cdot (y-v)|< |c_v|/10
$. It implies that $|\{y\in U_v : |\gamma(y) + u\cdot y - \beta| \le \lambda^{-1}\}|=0$ whenever $|\widetilde \beta|\ge|c_v|/5$ and $\lambda\ge 100|c_v|^{-1}$. Since \eqref{ineq:sublevelset} is trivial when $\lambda\le 100|c_v|^{-1}$, to establish the desired inequality, it is enough to consider the case $$|\widetilde \beta|<|c_v|/5.$$

Under this assumption, we obtain $|c_v - \widetilde u\cdot v - \widetilde \beta|\sim |c_v|$, which ensures that the argument in the proof of Theorem \ref{thm_lin_pert} remains valid with the modified choice
\[
    t=-x_d/(c_v - \widetilde u\cdot v - \widetilde \beta).
\] 
 Note that for this $t$, we have
\[
x_d -t(\gamma_v(y) + \widetilde u\cdot (y+ v) -c_{v})=-t(\gamma_v(y) + \widetilde u\cdot y - \widetilde \beta).
\] 
With the choice of $t$, for the function $f$ given by \eqref{def:ff}, we can check that  $\widetilde M_{\widetilde u,v} f$ is greater than
    \[
        \mathds1_{\mathbb B_{1/2}^{d-1}(0)}(x') \mathds1_{[-10^{-1}, 0)}(x_d/c_v) \int_{U_v'} (t|\gamma_v(y) + \widetilde u\cdot y - \widetilde \beta|)^{-1/p} \,\mathrm dy.
    \]
    Indeed, for $x$ satisfying  $x_d/c_v\in [-10^{-1},0)$, one has \eqref{ineq:x'} and
    \begin{align*}
        \mathds1_{[-1, 1]}(-t (\gamma_v(y)+ { \widetilde u \cdot y - \widetilde \beta})/c_v) \ge 1,
    \end{align*}
    whenever $\widetilde u\in\mathbb B_{|c_v|/10}^{d-1}(0)$, $|\widetilde \beta|\le 5^{-1}|c_v|$, and $y\in U_v'$.
    Then, by the $L^{p_0}$ boundedness of $M_\Gamma[\mathds1_U]$, it follows that
    \begin{align}
    \nonumber
        \int_{U_v'} | \gamma_v(y)+\widetilde u\cdot y - \widetilde \beta|^{-1/p}  \,\mathrm dy \le C |c_v|^{-\frac1p}.
    \end{align}
    The desired sublevel set estimate follows by applying Chebyshev's inequality combined with \eqref{identity:gammaub}.
\end{proof}

\subsection{Stationary set method and proof of Theorem~\ref{thm_Fdecay2}}

In this section, we introduce a method developed by Basu--Guo--Zhang--Zorin-Kranich \cite{BGZZ}, which helps us to control the following oscillatory integral,
\begin{equation}
    \label{oscillatory}
    \int_{[-1,1]^n} e^{i\lambda\phi(x)}\,\mathrm dx
\end{equation}
in terms of the size of a \emph{mid-level set},
\begin{equation}
\label{sublevelset}
|\{x\in [-1,1]^n: |\phi(x) - a|\le 1/\lambda\}|.   
\end{equation}
Then, together with Theorem~\ref{thm_lin_pert2}, we prove Theorem~\ref{thm_Fdecay2}.

In \cite{BGZZ}, the case where $\phi$ is semi-algebraic was considered in order to establish uniform bounds with respect to $\phi$. In contrast, we are concerned with linear perturbations of the form $\phi(x)+w\cdot x$, and our goal is to obtain uniform bounds with respect to $w$ only, rather than $\phi$. Since we do not pursue uniform estimates with respect to $\phi$, the estimate in \cite{BGZZ} can be extended to general analytic functions $\phi$. More precisely, the following holds.
\begin{theorem}\label{thm:BGZZ}
    Let $\phi: [-2,2]^n\rightarrow\mathbb R$ be an analytic function.  Then, there is a constant $C>0$ such that for all $\lambda\ge1$ and $w\in \bR^n$,
    \[
    \Big|\int_{[-1,1]^{n}} \mathrm{e}^{i\lambda(\phi(x)+w\cdot x)}~\mathrm{d}x \Big|\le C \sup_{a\in\bR}|\{x\in [-1,1]^{n}:|\phi(x)+w\cdot x-a|\le 1/\lambda\}|.
    \]
\end{theorem}

Theorem \ref{thm:BGZZ} can be deduced from \cite[Theorem 3.7]{LO}, where the last author extended the results in \cite{BGZZ} to the restricted analytic setting. In particular, it was proved that \eqref{oscillatory} can be bounded by \eqref{sublevelset} when $\phi$ is analytic with definable perturbations. Thus, Theorem \ref{thm:BGZZ} is a direct consequence of \cite[Theorem 3.7]{LO}. For the reader's convenience, however, we provide a self-contained proof below.

Before proving Theorem \ref{thm:BGZZ}, we first show the following lemma, which is a crucial observation for the argument.
\begin{lemma}\label{lem_monotonicity}

  Let $\phi$ be an analytic function defined on $[-2,2]^n$ and $a, \delta \in \mathbb R$.
  Let $W$ be given by
  \[W:=\{(x,w,a,\delta)\in[-1,1]^n\times \bR^n\times  \bR^2:|\phi(x)+w\cdot x-a|\le \delta \}.\]
   We define 
    \[        A(a,w,\delta) = \int_{\bR^n} \mathds1_{W}(x,w, a,\delta)~\mathrm{d}x,    \]
    and $N(w, \delta)$ is the number of times  $a\mapsto A(a,w,\delta)$ changes monotonicity. Then $\sup_{w,\delta} N(w,\delta)<\infty$.
\end{lemma}
\begin{proof}
    The lemma is a consequence of Proposition \ref{prop:BGZZ}. Indeed, we can rewrite $\mathds 1_W$ as follow:
    \[
    \mathds1_W= \mathds1_{[-1,1]^n}(x)\prod_{*=\pm} \mathds1_{\{f_{*}\ge0\}},
    \]
    where 
    \[
f_{\pm}=\pm(\phi(x) \mathds1_{[-1,1]^{n}}+w\cdot x-a)+\delta.    \]
 Thus, $\mathds1_W$ is a finite product of definable functions in $\bR_{an}$, which implies $\mathds1_W$ is definable in $\bR_{an}$. By Proposition \ref{prop:LO}, it follows that $A$ is definable in $\bR_{an,exp}$. Therefore, Proposition \ref{prop:BGZZ} gives the desired.
\end{proof}

Now, we prove Theorem \ref{thm:BGZZ} using Lemma \ref{lem_monotonicity}.
\begin{proof}[Proof of Theorem \ref{thm:BGZZ}]
    We follow the proof of Basu--Guo--Zhang--Zorin-Kranich \cite{BGZZ}.
    Let $\phi_w^\lambda(x):=\lambda(\phi(x)+w\cdot x)$.
    Then we have
    \[
        \Big| \int_{-1}^{1} \mathrm{e}^{i a}~\mathrm{d}a \times \int_{[-1,1]^n} \mathrm{e}^{i \phi_w^\lambda(x)} ~\mathrm{d}x\Big|
        \sim \Big| \int_{[-1,1]^n} \mathrm{e}^{i \phi_w^\lambda(x)} ~\mathrm{d}x\Big|.
    \]
    By Fubini's theorem, it follows that
    \begin{align*}
        \Big| \int_{[-1,1]^n} \mathrm{e}^{i \phi_w^\lambda(x)} ~\mathrm{d}x\Big|
        \sim &\Big|\int_{[-1,1]^n} \int_{\phi_w^\lambda(x)-1}^{\phi_w^\lambda(x)+ 1} \mathrm{e}^{i a}~\mathrm{d}a\mathrm{d}x \Big|\\
        =&\int_{\bR} |\{x\in [-1,1]^n : |\phi_w^\lambda(x)-a|\leq 1\}| \mathrm{e}^{ia}~\mathrm{d}a.
    \end{align*}
    Let $S_{\phi}(a;w, \lambda)$ denote the size of a stationary set of $\phi^\lambda_w$ of height $a$ and thickness $1$,
    \[
        S_\phi(a;w, \lambda) = |\{x\in [-1,1]^n : |\phi_w^\lambda(x)-a|\leq 1\}|.
    \]
    By Lemma~\ref{lem_monotonicity} with $A(a/\lambda,w,\lambda^{-1})=S_\phi(a;w,\lambda)$, we know that $S_\phi(a;w,\lambda)$ changes its monotonicity only $O_{d, \phi}(1)$ times in terms of $a$.
    Let $I_j$'s be disjoint intervals on which $S_\phi(\cdot;w,\lambda)$ is monotonic.
    Then, we show
    \begin{align*}
        \Big| \int_{\bR}  S_\phi(a; w,\lambda)  \mathrm{e}^{ia}~\mathrm{d}a \Big|
        = &\sum_j  \Big| \int_{ I_j} S_\phi(a;w, \lambda) \mathrm{e}^{ia}~\mathrm{d}a \Big|\\
        \le & C \sup_{a\in\bR} S_\phi(a; w,\lambda).
    \end{align*}
    Indeed, the last inequality follows from integration by parts, as the total variation of $S_\phi(a; w,\lambda)$ on each $I_j$ is bounded by $2\sup_{a\in\bR}S_\phi(a; w,\lambda)$. 
    Hence, we have shown that
    \[
        \Big| \int_{[-1,1]^n} \mathrm{e}^{i \lambda(\phi(x)+w\cdot x)} ~\mathrm{d}x\Big| \le C \sup_{a\in\bR} S_\phi(a; w,\lambda).
    \]
    The desired result follows by substituting $a$ by $\lambda a$, since $S_\phi(\lambda a; w,\lambda)=|\{x\in [-1,1]^n:|\phi(x)+w\cdot x -a|\le 1/\lambda]\}|$.
\end{proof}

Now, we prove Theorem~\ref{thm_Fdecay2}.
\begin{proof}[Proof of Theorem~\ref{thm_Fdecay2}]
     Recall that $\mathcal F [\mu[\psi]](\xi', \xi_d)$ is given by
\[
    \int_{\mathbb R^{d-1}} \mathrm e^{-2\pi i (\xi', \xi_{d})\cdot (y, \gamma(y))}\psi(y)\,\mathrm dy.
\]
Let $|(\xi', \xi_d)| \sim \lambda$.
If one has $|\xi'| \ge C |\xi_d|$ with $C$ depending on $|\nabla \gamma|$, then the desired decay is obtained via integration by parts.
Thus, we assume that $|\xi'| \lesssim |\xi_d| \sim \lambda$.
Without loss of generality, we take $\lambda = \xi_d$ so that the matter is reduced to proving the following oscillatory integral  estimate:
\[
    \Big|\int \mathrm{e}^{i\lambda(\gamma(x)+u\cdot y)}\psi(y)~\mathrm{d}y\Big|\lesssim (c_\psi \lambda)^{-1/p},
\]
where $u= \lambda^{-1}\xi'$. 

By Theorem \ref{thm_lin_pert2} with $p>p_0>p_{cr}$, for each $v\in \supp (\psi)$ there exists a cube of side length $r_v>0$, $$U_v=v+r_v[-1,1]^{d-1},$$ satisfying \eqref{ineq:sublevelset}. Then, we choose a finite subcover $\{U_{v_i}\}_{i=1}^M$ for $v_i\in \supp(\psi)$ such that $\supp(\psi)   \subset \cup_{i=1}^M U_{v_i}$. Using a partition of unity $\{\varphi_{v_i}\}_{i=1}^M$ subordinate to $\{U_{v_i}\}_{i=1}^M$, it suffices to show that
\[
    \Big|\int \mathrm{e}^{i\lambda(\gamma(x)+u\cdot y)}\psi(y)\varphi_{v_i}(y)~\mathrm{d}y\Big|\lesssim (c_\psi \lambda)^{-1/p}.
\]
Using Fourier series of $\psi\varphi_{v_i}$, one has
\[
    \psi\varphi_{v_i}(y) = \sum_{m\in\mathbb Z^{d-1}} c_n \mathrm e^{i m\cdot y},\quad\text{with}\quad  \sum_m |c_m| \lesssim 1.
\]
Then we obtain
\[
    \Big|\int \mathrm{e}^{i\lambda(\gamma(y)+u\cdot y)}\psi(y)\varphi_{v_i}(y)\,\mathrm{d}y\Big|
     \lesssim \sup_{u\in\mathbb R^{d-1}} \Big|\int_{  U_{v_i}} \mathrm{e}^{i\lambda(\gamma(y)+u\cdot y)}~\mathrm{d}y\Big|.
\]
By \eqref{ineq:sublevelset} and the lower bound $c_\psi\le |c_{v_i}|$, it is enough to show that for all $i=1,\cdots, M$, 
\begin{equation}
\label{ineq:subleveluvi}
        |\int_{  U_{v_i}} \mathrm{e}^{i\lambda(\gamma(y)+u\cdot y)}~\mathrm{d}y\Big|\lesssim \sup_{a}|\{y\in  U_{v_i}:|\gamma(y)+u\cdot y-a|\le 1/\lambda\}|.
\end{equation}
To verify this for fixed $i$, we use change of variables, $y\mapsto r_{v_i}y+v_i$. This yields that
\[
|\int_{  U_{v_i}} \mathrm{e}^{i\lambda(\gamma(y)+u\cdot y)}~\mathrm{d}y|=r_{v_i}^{d-1}|\int_{[-1,1]^{d-1}}e^{i\lambda(\gamma(v_i+r_{v_i}y)+u\cdot (r_{v_i}y))} \mathrm{d}y|.
\]
Applying Theorem \ref{thm:BGZZ} with 
\[
\phi(y)=\gamma(v+r_{v_i}y), \quad w=r_{v_i} u,
\]
 we obtain  that
\[
\begin{aligned}
    |\int_{  U_{v_i}} \mathrm{e}^{i\lambda(\gamma(y)+u\cdot y)}~\mathrm{d}y|\lesssim r_{v_i}^{d-1}\sup_{a}|\{y\in  [-1,1]^{d-1}:|\phi(y)+w\cdot y-a|\le 1/\lambda\}|.
\end{aligned}
\]
Finally, using the scaling identity $|rE+v|=r^{d-1}|E|$ for a set  $E\subset \bR^{d-1}$, we conclude that the desired estimate \eqref{ineq:subleveluvi} holds.
\end{proof}

\appendix

\section{Non-transversality with scaling structure}\label{sec_scaling}

The purpose of this appendix is to provide short proofs of results analogous to Theorem~\ref{thm_main} where the hypersurfaces enjoy certain scaling structures, specifically when $\Gamma\subset \bR^2$ is a smooth finite type curve or when $\Gamma$ is a smooth convex finite line type hypersurface. 
If a hypersurface $\Gamma$ has an inherent scaling structure, an analogue of Theorem~\ref{thm_main} remains valid even without analyticity.
We first consider the case in $\mathbb R^2$.
\begin{theorem}
    Let $\gamma:\bR\rightarrow \bR$ be a smooth function of finite type and $(y_{nt},\gamma(y_{nt}))$ be a non-transversal point of $\Gamma$. If $\psi$ has a sufficiently small support near $y_{nt}$, then $M_\Gamma[\psi]$ is bounded on $L^p(\bR^2)$ for all $p>2$. 
\end{theorem}
\begin{proof}
After a suitable rotation, we may assume that $\gamma(y_{nt})=0$ and $\nabla\gamma(y_{nt})=0$ and that $\supp(\psi)\subset [y_{nt}-1,y_{nt}+1]$. Since $\gamma$ is of finite type, there exists $k\ge2$ such that
\[
\gamma(y)=c_\gamma(y-y_{nt})^k+R(y-y_{nt})
\]
for some $c_\gamma$, where the remainder $R$ satisfies that $|(d/dy)^m R(y)|\lesssim  |y|^{k+1-m}$ for $0\le m\le k+1$ and $|y|\le 1$. By a scaling, one can assume that $c_\gamma=1$.

Define
\[
\psi_\ell(y)=\psi(2^{-\ell}y+y_{nt})\beta_1(y).
\] 
Since $\supp(\psi)\subset [y_{nt}-1,y_{nt}+1]$, we have the decomposition $\psi(y)=\sum_{\ell\ge0}\psi_\ell(2^\ell(y-y_{nt}))$. This decomposition and change of variables $y\mapsto 2^{-\ell}y+y_{nt}$ yield that
    \[
    \begin{aligned}
    A_\Gamma&[\psi]f(t,x_1,x_2)\\
    &=\int f(x_1-ty,x_2-t\gamma(y))\psi(y)dy     \\
    &=\sum_{\ell=0}^\infty \int f(x_1-ty,x_2-t\gamma(y))\psi_\ell(2^\ell(y-y_{nt}))dy \\
    &=\sum_{\ell=0}^\infty 2^{-\ell}\int f(x_1-t(2^{-\ell}y+y_{nt}),x_2-t(2^{-\ell k}y^k+R(2^{-\ell}y)))\psi_\ell(y)dy.
    \end{aligned}
    \]
    By changing of variables $(x_1,x_2)\mapsto (2^{-\ell}x_1,2^{-\ell k}x_2)$, one can observe that
    \[
    \|M_\Gamma[\psi]\|_{p\rightarrow p}\le \sum_{\ell=0}^\infty   2^{-\ell}  \|\widetilde {M_\ell}\|_{p\rightarrow p},
    \]
    where
    \[
    \widetilde {M_\ell}f(t,x)=\sup_{t>0}\big|\int _{|y|\sim 1}f(x_1-t(y+2^\ell y_{nt}),x_2-t(y^k+2^{\ell k}R(2^{-\ell}y)))\psi_\ell(y)dy\big|.
    \]

By the property of $R$, there exists a positive integer $L$ such that
\[
2^{\ell(k-2)}R''(2^{-\ell}y)\le \frac{1}{100},
\]
for all $|y|\le 2$ and all $\ell\ge L$. Now, assume that $\supp(\psi)$ is contained in a small neighborhood $[y_{nt}-2^{-L},y_{nt}+2^{-L}]$. Then $\widetilde{M_\ell}=0$ for $\ell<L$, so it is enough to show that 
\[
\sum_{\ell=L}^\infty   2^{-\ell}  \|\widetilde {M_\ell}\|_{p\rightarrow p}\lesssim 1.
\]
    For $\ell\ge L$, we have a favorable lower bound on the second derivative of the function $y^k+2^{\ell k}R(2^{-\ell}y)$. Due to this observation, we can apply
     the following estimate (\cite[(16)]{IS2}, modified version of Sogge's theorem \cite{Sogge}), 
    \[
    \|M_\Gamma f\|_p\lesssim d(0,\Gamma)^{1/p}\|f\|_p \quad \text{ for } p>2.
    \]
    Since $\widetilde {M_\ell}=M_{\Gamma_\ell}[\psi_\ell]$ where $\Gamma_\ell=\{(y+2^{\ell}y_{nt},y^k+2^{\ell k}R(2^{-\ell}y)):|y|\sim 1\}$ and $\psi_\ell$ satisfies a uniform bound on its derivatives, Sogge's theorem implies that
\[
\|\widetilde {M_\ell} \|_{p\rightarrow p}\lesssim 2^{\ell/p}, \quad \text{for }p>2.
\]
 Since $p>2$, it deduces $\sum_{\ell=L}^\infty 2^{-l}\Vert \widetilde M_\ell\Vert_{p\to p}\lesssim 1$ and we obtain the desired estimate.
\end{proof}

We also obtain analogous results for convex hypersurfaces in $\mathbb R^d$.
\begin{theorem}
    Let $\gamma:\bR^{d-1}\rightarrow \bR$ be a smooth convex function of finite line type and $(y_{nt},\gamma(y_{nt}))$ is a non-transversal point of $\Gamma$. If $\psi$ has a sufficiently small support near $y_{nt}$, then $M_\Gamma[\psi]$ is bounded on $L^p(\bR^d)$ for all $p>2$. 
\end{theorem}
Here, we can deal with smooth hypersurfaces as well as analytic surfaces. The proof below follows the arguments in \cite{IS2}.
\begin{proof}
As in the case $d=2$, we first assume that $\gamma(y_{nt})=0$ and $\nabla\gamma(y_{nt})=0$. 
By the result of Schulz \cite{schulz1991}, there exist a  nontrivial multi-homogeneous function $Q$ and a remainder term $R$ such that
\[
\gamma(y+y_{nt})=Q(y)+R(y).
\]
More precisely, there exist integers $k_1,\cdots, k_{d-1}\ge2$ such that
 $$Q(T_\ell y)=2^{-\ell}Q(y),$$ and $2^{\ell}R(T_\ell y)\ll1$ for sufficiently large $\ell$ and $|y|\le 2$ where the anisotropic dilation $T_\ell$ is defined by $$T_\ell y=(2^{-\ell/k_1}y_1,\cdots,2^{-\ell/k_{d-1}}y_{d-1}).$$

Using dyadic decomposition, for $(\tilde x,x_d)\in \bR^{d-1}\times \bR$, we get an analogous identity,
    \[
    \begin{aligned}
    &A_\Gamma[\psi]f(t,\tilde x,x_d)\\
    &=\sum_{\ell=1}^\infty 2^{-\ell(\sum_{i=1}^{d-1}\frac 1{k_i}
    )}\int f(\tilde x-t(T_\ell y+y_{nt}), x_d-t(2^{-\ell}Q(y)+R(T_\ell y)))\psi_\ell(y)dy,
    \end{aligned}
    \]
    where 
$$\psi_\ell(y)=\psi(T_\ell y+y_{nt})\beta_*(y)$$
and $\beta_*$ is a smooth function supported in the set $[-2,2]^{d-1}\setminus [-1,1]^{d-1}$.
    By a change of variables $(\tilde x,x_d)\mapsto (T_\ell x,2^{-\ell }x_d)$, it suffices to show that there exists a positive integer $L$ such that
    \[
\sum_{\ell=L}^\infty   2^{-\ell(\sum_{i=1}^{d-1}\frac 1{k_i}
    )}  \|\widetilde {M_\ell}\|_{p\rightarrow p}\lesssim 1,
    \]
    where
    \[
    \widetilde {M_\ell}f(x)=\sup_{t>0}|\int f(x-t(y+T_{-\ell}y_{nt}),x_d-t(Q(y)+2^{\ell}R(T_\ell y)))\psi_\ell(y)dy|.
    \]
As proved in \cite[Section 3]{IS2},
the graph $(y+T_{-\ell}y_{nt},Q(y)+2^{\ell}R(T_\ell y))$ has at least one nonvanishing principal curvature on $\supp(\psi_\ell)$, if $\ell$ is sufficiently large. Thus, Sogge's theorem implies
\[
\|\widetilde {M_\ell}\|_{p\rightarrow p}\lesssim \max_i 2^{\frac{\ell}{pk_i}}.
\]
It gives the desired bound, since
\[
\frac{1}{p\min_i k_i}-\sum_{i=1}^{d-1}\frac 1{k_i}<0
\]
for every $p>1$.
\end{proof}

\section*{Acknowledgement}
The authors would like to thank Alex Iosevich for valuable comments.
This work is supported by the National Research Foundation of Korea (RS-2021-NR061906; J. B. Lee), (RS-2024-00461749; J. B. Lee and J. Oh); KIAS Individual Grant (MG098901; J. Lee), (SP089101; S. Oh).

\end{document}